\documentclass[pdflatex,sn-mathphys-num]{sn-jnl}
\usepackage{graphicx}%
\usepackage{multirow}%
\usepackage{amsmath,amssymb,amsfonts}%
\usepackage{amsthm}%
\usepackage{mathrsfs}%
\usepackage[title]{appendix}%
\usepackage{xcolor}%
\usepackage{textcomp}%
\usepackage{manyfoot}%
\usepackage{booktabs}%
\usepackage{algorithm}%
\usepackage{algorithmicx}%
\usepackage{algpseudocode}%
\usepackage{listings}%
\usepackage{epstopdf}
\usepackage{mathtools}
\usepackage{setspace} 
\usepackage{url}
\usepackage{subcaption}
\usepackage{comment}
\usepackage{soul}
\let\orcidlogo\undefined 
\usepackage{orcidlink}

\newcommand{\SUT}{\operatorname{SUT}}
\newcommand{\bu}{\mathbf{u}}   
\newcommand{\bhB}{\widehat{B}}
\newcommand{\bhU}{\widehat{U}}
\newcommand{\hu}{\widehat{u}}
\newcommand{\bhV}{\widehat{V}}
\newcommand{\hv}{\widehat{v}}
\newcommand{\hx}{\widehat{x}}

\newcommand{\bub}{\delta b}
\newcommand{\halpha}{\widehat{\alpha}}
\newcommand{\hbeta}{\widehat{\beta}}
\theoremstyle{thmstyleone}
\newtheorem{theorem}{Theorem}[section]
\newtheorem{corollary}{Corollary}[section]
\newtheorem{lemma}{Lemma}[section]
\newcommand\Tstrut{\rule{0pt}{2.6ex}}       
\renewcommand{\tabcolsep}{2mm} 
\newcommand{\RO}{reorthogonalisation }
\newcommand{\norm}[1]{\left\|#1\right\|}

\begin{document}
\title[Reorthogonalisation in LSQR]
{Revisiting the role of reorthogonalisation in LSQR}

\author[1]{\fnm{Ieva} \sur{Dau\v{z}ickait\.{e}}} \email{ieva.dauzickaite@cerfacs.fr \orcidlink{0000-0002-1285-1764}}
\author[2]{\fnm{Melina A.} \sur{Freitag}}\email{melina.freitag@uni-potsdam.de \orcidlink{0000-0002-4539-2162}}
\author[1]{\fnm{Selime} \sur{G\"{u}rol}}\email{selime.gurol@cerfacs.fr \orcidlink{0000-0003-1269-4022}}
\author[3,4]{\fnm{Amos S.} \sur{Lawless}}\email{a.s.lawless@reading.ac.uk \orcidlink{0000-0002-3016-6568}}
\author[5]{\fnm{Alison} \sur{Ramage}}\email{a.ramage@strath.ac.uk \orcidlink{0000-0003-4709-0691}}
\author*[3,6]{\fnm{Jennifer A.} \sur{Scott}}\email{jennifer.scott@reading.ac.uk \orcidlink{0000-0003-2130-1091}}
\author[7]{\fnm{Jemima M.} \sur{Tabeart}}\email{j.m.tabeart@tue.nl \orcidlink{0000-0001-6806-8608}}

\affil[1]{\orgname{CERFACS},
\orgaddress{\city{Toulouse}, \country{France}}}
\affil[2]{\orgdiv{Institute for Mathematics}, \orgname{University of Potsdam}, 
\orgaddress{\city{Potsdam}, \country{Germany}}}
\affil[3]{\orgdiv{School of Mathematical, Physical and Computational Sciences}, \orgname{University of Reading}, 
\orgaddress{\city{Reading}, \country{England}}}
\affil[4]{\orgdiv{National Centre for Earth Observation}, \orgname{University of Reading},
\orgaddress{\city{Reading}, \country{England}}}
\affil[5]{\orgdiv{Department of Mathematics and Statistics}, \orgname{University of Strathclyde}, 
\orgaddress{\city{Glasgow}, \country{Scotland}}}
\affil[6]{\orgdiv{Scientific Computing Department}, \orgname{Rutherford Appleton Laboratory, Science \& Technology 
Facilities Council}, 
\orgaddress{\city{Didcot}, \country{England}}}
\affil[7]{\orgdiv{Department of Mathematics and Computer Science}, \orgname{Eindhoven University of Technology},
\orgaddress{\city{Eindhoven}, \country{The Netherlands}}}

\abstract{
In this paper, we revisit the role of reorthogonalisation in the LSQR
method for large-scale overdetermined linear least-squares problems. Motivated
primarily by practical applications where matrix--vector products dominate the
computational cost and only a limited number of iterations can be afforded because of time and/or memory restrictions, we
investigate when reorthogonalisation can significantly accelerate convergence. 
Through synthetic test problems with controlled singular value distributions and
experiments on practical sparse problems, we demonstrate that reorthogonalisation is
particularly beneficial when the least-squares matrix has large
outlying singular values.  In such cases, loss of orthogonality can
severely delay convergence. We
explore the effectiveness of full (one-sided) \RO when using LSQR,
and present alternative  \RO strategies that can provide substantial
acceleration compared to no \RO and significant storage savings compared to full reorthogonalisation.
The possibility of using low precision arithmetic for the \RO is also 
investigated. We {present} bounds that demonstrate the quality of the solution depends on the level of orthogonality maintained by the selected \RO strategy. Our results suggest that reorthogonalisation should be
regarded as a practical tool for modern least-squares applications that can be employed 
instead of, or in addition to,  preconditioning. }

\keywords{Linear least-squares problems, sparse problems, Krylov subspace methods,
reorthogonalisation, mixed precision, LSQR. \\ 
MSC codes 65F10, 65F20, 65F25 }


\maketitle

\section{Introduction}
\subsection{Motivation and context}
Large-scale linear least-squares problems arise throughout scientific computing, including in inverse problems, data assimilation, machine learning, and
PDE-constrained optimisation. They lead to overdetermined systems of the form
\begin{equation}
\min_{x \in \mathbb{R}^n} \| b- Jx \|_2,
\label{eq:LSproblem}
\end{equation}
where $J \in \mathbb{R}^{p \times n}$ with $p \gg n$. In practical applications, $J$ may not be available explicitly and matrix--vector 
products with $J$ and $J^T$ are often significantly more expensive than vector 
operations. In such settings, Krylov subspace methods, notably
CGLS~\cite{hest:52}, LSQR~\cite{pasa:82} and LSMR~\cite{fosa:2011}, are attractive because of 
their matrix-free structure, modest storage requirements, and suitability for large problems.

In exact arithmetic, Krylov subspace methods generate
mutually orthonormal  basis vectors. In finite precision arithmetic, 
however, orthogonality is  lost, leading to a slowdown in convergence (or even a failure to achieve the requested accuracy). The loss of orthogonality can happen gradually or relatively rapidly. While this phenomenon is well known and has been 
extensively analysed in the context of eigenvalue computations \cite{ paige1971computation, ParlettScott1979} and GMRES-type 
methods \cite{saad2003,saad1986}, its impact on least-squares solvers is often regarded as 
secondary. Consequently, many of the readily available implementations 
of Krylov methods for least-squares problems (including the LSQR packages offered in Matlab by MathWorks\footnote{\url{https://uk.mathworks.com/help/matlab/ref/lsqr.html}}, in Python by SciPy\footnote{\url{https://docs.scipy.org/doc/scipy/reference/generated/scipy.sparse.linalg.lsqr.html}}, and
in a number of languages in Stanford University's Systems Optimisation Library\footnote{\url{https://web.stanford.edu/group/SOL/software/lsqr/}}) do not include options for
reorthogonalisation; this may be why many 
of those in the numerical linear algebra
community do not routinely employ it. 

Modern large-scale applications challenge this conventional view. In many 
contemporary inverse and data assimilation problems, matrix--vector products 
with $J$ and $J^T$ dominate the computational cost, preconditioning options are few 
or unavailable, and only a restricted number of iterations can be afforded. In 
some instances, high solution accuracy may be required. In such regimes, slow 
convergence of the iterative method can cause substantially more expense than additional vector 
orthogonalisation. Furthermore, the spectral properties of $J$ in these 
applications often include large outlying singular values
that can degrade the effectiveness of 
Krylov subspace approximations when orthogonality is lost, delaying convergence.
Although preconditioning may speed up convergence, so that \RO
is typically considered to be less crucial, for some problems preconditioning does not do enough, or may not be
readily available. 

This paper revisits the role of reorthogonalisation in least-squares Krylov subspace methods from a modern computational perspective, focussing specifically on 
the widely-used LSQR method. When iterative methods are applied to the normal equations, finite-precision effects are especially pronounced: the condition number of $J^TJ$ is the square of that of $J$, which amplifies the impact of rounding errors. The challenge of problems with large condition numbers has helped motivate extensive research into reorthogonalisation techniques for solvers for linear systems of equations, including the Conjugate Gradient (CG) method~\cite{meurant_lanczos_2006, Notay1993}.
CG applied to the normal equations and LSQR are mathematically equivalent in exact arithmetic~\cite{pasa:82} but their finite-precision behaviour can differ considerably. Although 
for least-squares problems LSQR has been recognised as numerically superior to CG for many years, the numerical properties of LSQR have received comparatively less attention. We seek to investigate the circumstances under which loss of orthogonality significantly degrades its convergence, with particular attention to regimes where CG is known to require reorthogonalisation~\cite{carson2024}. We examine how different partial reorthogonalisation strategies can reduce iteration counts, memory requirements, and overall computational cost, and explore how the effectiveness of these strategies relates to the distribution of the singular values.
We also present forward-backward error bounds for the Golub--Kahan--Lanczos bidiagonalisation and a bound on the LSQR relative solution error that apply for a wide range of one-sided reorthogonalisation strategies.

Through a combination of synthetic test problems
with controlled spectral properties and benchmark matrices, we demonstrate that 
reorthogonalisation can substantially accelerate LSQR convergence when $J$ has 
large outlying singular values. In such cases, 
LSQR is sensitive to 
loss of orthogonality, which can result in convergence delay in finite
precision computations. We
further show that full (one-sided) reorthogonalisation is often unnecessary: appropriately 
chosen partial strategies can recover most of the acceleration at considerably 
reduced cost (in terms of time and storage).
These findings suggest that reorthogonalisation should not be viewed solely as 
a theoretical safeguard for maintaining orthogonality, but rather as a practical 
mechanism for accelerating convergence in large-scale least-squares 
applications, and may be particularly useful when effective preconditioners are unknown and/or matrix--vector products are very
expensive, limiting the number of iterations that can be afforded.

The remainder of the paper is organised as follows. In Section \ref{sec:lsqr}, we state the LSQR algorithm and 
summarise some reorthogonalisation strategies used in the literature, alongside new error bounds for the Golub--Kahan--Lanczos bidiagonalisation. 
Section \ref{sec:Ritz} discusses the role of singular values and Ritz values in LSQR convergence using 
synthetic, but insightful, examples. Techniques for implementing reorthogonalisation in practical settings, where only 
a fixed memory budget is available, are discussed in Section \ref{sec:local}. Section \ref{sec:precision} 
provides a framework for mixed precision computation within reorthogonalisation in the LSQR
method, including a theoretical bound on convergence, before we 
finish with some concluding remarks.

\subsection{ Notation and test environment}
The following is a summary of the key notation used throughout 
the paper. We write $\|\cdot\|=\|\cdot\|_2$ to denote the 
Euclidean norm, and $\kappa(A)=\kappa_2(A)=\|A\|_2\|A^\dagger\|_2$ denotes the 
2-norm condition number of a general matrix $A$ where $A^\dagger$ is the Moore--Penrose pseudoinverse of $A$. 
Subscripts denote Krylov iteration 
numbers, so that $x_j$ is the approximation to the true solution 
$x_{true}$ at iteration $j$. For a square matrix $A$ and vector $y$, the Krylov subspace of order $j$ is defined by
\[
\mathcal{K}_j(A,y)=\operatorname{span}[\,y,Ay,\ldots,A^{j-1}y\,].
\]
We write $A_j=[a_1,\ldots,a_j]$ for a matrix with $j$ columns, and define the orthogonality level of its columns using
\begin{equation}\label{eq:orth index}
\mu(A_j)=\norm{\SUT\!\left(I_j-A_j^TA_j\right)}.
\end{equation}
Here $\SUT(\cdot)$ denotes the strictly upper triangular part of a matrix and $I_j$ is the identity matrix of order $j$ with columns $e_i^{(j)}$,  $1\le i\le j$. Finally, $\mathcal{S}^A$ denotes an index set corresponding to a subset of the columns of $A$.

Unless stated otherwise, our experiments are conducted in MATLAB version 2024b.
 LSQR with modified Gram-Schmidt reorthogonalisation was implemented following \cite{benb:1999}, with the stopping criteria 
(see Section~\ref{stopping}) using the open source MATLAB code associated with \cite{pati:2024}.

\section{LSQR with reorthogonalisation}
\label{sec:lsqr}
The Krylov subspace method 
LSQR for solving (\ref{eq:LSproblem}) is based on the Golub--Kahan--Lanczos 
bidiagonalisation\footnote{Sometimes it is referred to as Golub-Kahan bidiagonalisation, or simply as Lanczos bidiagonalisation.} (GKLB) of
$J$; {see the original paper \cite{pasa:82} and the book \cite{Bjorckbook2024} for an overview and further references}. Let $r_0 = b - Jx_0$ be the initial residual.
In exact arithmetic, at iteration $j$ the process generates orthonormal vectors
$u_{j+1}$, $v_{j+1}$ and scalars $\alpha_{j+1}$, $\beta_{j+1}$ such that
\[
J v_j = \alpha_j u_j + \beta_{j+1} u_{j+1}, \qquad
J^T u_{j+1} = \beta_{j+1} v_j + \alpha_{j+1} v_{j+1}.
\]
\begin{algorithm}
\caption{LSQR with optional reorthogonalisation}
\label{alg:LSQR}
Input: $J$, $b$, initial guess $x_0$, initial index sets $\mathcal{S}^U$ and $\mathcal{S}^V$ (which may be empty sets).\\
Output: Approximate least-squares solution $x$.
\vskip+2mm\hrule\vskip+2mm
\setstretch{1.35}
\begin{algorithmic}[1]
\State $r_0 = b - Jx_0$\Comment{Initialisation}
  \State $\beta_{1} = \norm{r_0}$, \quad $u_{1}= r_0/ \beta_{1}$
\State $q = J^T u_1$, \quad $\alpha_{1} = \norm{q}$, \quad $v_{1}= q/ \alpha_{1}$
\State $w_1 = v_1$, \quad $\bar\phi_1 = \beta_1$, \quad $\bar\rho_1 = \alpha_1$
\For{$j = 1,2,\dots$ until convergence}
  \State $q = J v_j - \alpha_j u_j$ \label{step:q_for_u}
  \State Update $\mathcal{S}^U$ and reorthogonalise $q$ against each $u_i$ for which $i \le j$ and $i \in \mathcal{S}^U$ \Comment{Skip  if $\mathcal{S}^U = \emptyset$} \label{step:Ro_u}
  \State $\beta_{j+1} = \norm{q}$, $u_{j+1}= q/ \beta_{j+1}$ \label{step:get_beta}
  \State $q = J^T u_{j+1} - \beta_{j+1} v_j$ \label{step:q_for_v}
  \State Update $\mathcal{S}^V$ and reorthogonalise $q$ against each $v_i$ for which $i \le j$ and $i \in \mathcal{S}^V$ \Comment{Skip if $\mathcal{S}^V = \emptyset$} \label{step:Ro_v}
  \State  $\alpha_{j+1} = \norm{q}$, $v_{j+1}= q/ \alpha_{j+1}$ \label{step:get_alpha}
  \State $\rho_j = (\bar\rho_j^2 + \beta_{j+1}^2)^{1/2}$ \label{step:start_Givens}
  \State $c_j = \bar\rho_j / \rho_j$, \quad $s_j = \beta_{j+1} / \rho_j$ \Comment{Steps \ref{step:start_Givens} to \ref{step:end_Givens} are Givens rotations}
  \State $\theta_{j+1} = s_j \alpha_{j+1}$, \quad $\bar\rho_{j+1} = -c_j
\alpha_{j+1}$
  \State $\phi_j = c_j \bar\phi_j$, \quad $\bar\phi_{j+1} = s_j \bar\phi_j$ \label{step:end_Givens}
  \State $x_j = x_{j-1} + (\phi_j / \rho_j) w_j$                     \Comment{Solution update. If converged, return $x = x_j$}
  \State $w_{j+1} = v_{j+1} - (\theta_{j+1} / \rho_j) w_j$
\EndFor
\end{algorithmic}
\end{algorithm}

The columns of $V_j = [v_1,\dots,v_j]\in \mathbb{R}^{n\times j}$ and those of $U_j = [u_1,\dots,u_j]\in \mathbb{R}^{p\times j}$ span the Krylov subspaces
$\mathcal{K}_j(J^T J, J^T r_0)$ and 
$\mathcal{K}_j(JJ^T, r_0)$, respectively.  Defining 
the bidiagonal matrix 
\[B_j = \begin{bmatrix}
  \alpha_1 &          &       &             \\
  \beta_2 & \alpha_2  &       &             \\
          & \beta_3  & \ddots &             \\
          &         &  \ddots &   \alpha_j   \\
          &         &         & \beta_{j+1} \\  
\end{bmatrix}\in \mathbb{R}^{(j+1)\times j},\] 
the GKLB can be expressed via the matrix recurrences
\[
J V_j = U_{j+1} B_j, \qquad
J^T U_{j+1} = V_j B_j^T + \alpha_{j+1} v_{j+1} (e_{j+1}^{(j+1)})^T.
\]
GKLB is mathematically equivalent to  Lanczos tridiagonalisation applied to $C=\begin{bmatrix}
    0 &  J\\ J^T & 0
\end{bmatrix}$ with a starting vector 
$\begin{bmatrix}
    u_1^T & 0
\end{bmatrix}^T$ 
(see, for example, \cite{paige1974bidiagonalization,lars:98}). 
The $j$th iterate $x_j$ of LSQR can be expressed as a vector in the range of $V_j$ (possibly shifted by a nonzero initial guess)
that minimises the residual norm, that is,
\begin{displaymath}
    x_j = x_0+ \, V_j y_j, \quad \textrm{where} \quad 
    y_j =  \, \arg \min_y \norm{b- J (x_0+ V_j y)}
        =  \, \arg \min_y \norm{B_j y - \beta_1 e_1^{(j+1)}}. 
\end{displaymath}

The singular values of $B_j$ approximate the singular values of $J$.
These approximations are equal to the positive \emph{Ritz values} of $C$, that is, the optimal approximations of the positive eigenvalues of $C$ 
from the Krylov subspace generated by the Lanczos tridiagonalisation; see, e.g., \cite[Ch. 7]{Demmel1997}. The corresponding approximations 
to the eigenvectors are called \emph{Ritz vectors}.

In LSQR implementations, the approximate solution $x_j$ is updated recursively,
as outlined in Algorithm~\ref{alg:LSQR}. 
In this algorithm, 
$\mathcal{S}^U$ and $\mathcal{S}^V$ denote index sets that determine 
the choice of columns in $U$ and $V$ to be included in the reorthogonalisation at steps~\ref{step:Ro_u}
and \ref{step:Ro_v}. 
Choosing $\mathcal{S}^U$ and $\mathcal{S}^V$ to be empty sets corresponds to no reorthogonalisation.

\subsection{Stopping criteria for LSQR}
\label{stopping}
The choice of stopping criterion for LSQR is important and not straightforward.
The criteria proposed by Paige and Saunders~\cite{pasa:82} are widely used
but have well-documented limitations. In particular, the quantity
$\|J^T r_j\| / \|r_j\|$ can oscillate for
ill-conditioned problems and subsequently plateau, potentially failing to trigger
termination even when the iterate is already a backward stable
solution~\cite{cptp:2009}.
In our experiments, we use the stopping criterion recently proposed by
Pape\v{z} and Tich\'{y}~\cite{pati:2024}.
They show that a reliable estimate of
$\|x_{true} - x_j\|_{J^TJ}$, the energy norm of the solution error at the current iterate,
can be computed  at no additional cost from the scalars $\phi_j$ already produced
by the LSQR recurrence, and that this estimate remains
reliable in finite-precision arithmetic provided local orthogonality of
consecutive GKLB vectors is maintained and the LSQR solution has not 
reached its maximal attainable accuracy.
Following~\cite{sctu:2026}, we terminate LSQR at iteration $j$ if 
\begin{equation}
  \mathrm{ratio}_{PT} \equiv
  \frac{\mathrm{estim}_{\ell_j}}
       {\mathrm{estim}(\|J\|)\,\|x_j\| + \|b\|}
  \leq \delta,
  \label{eq:ratio_PT}
\end{equation}
where $\mathrm{estim}_{\ell_j}$ is the estimate of the energy norm
$\|x_{true} - x_{\ell_k}\|_{J^TJ}$ for an adaptively chosen recent iterate
$x_{\ell_k}$ where $k\le j$, and $\mathrm{estim}(\|J\|)$ is an estimate of $\|J\|$.
The normalisation makes $\mathrm{ratio}_{PT}$ dimensionless and ensures that
the tolerance $\delta$ has a consistent interpretation across problems of
different scales.
This criterion provides a direct and reliable measure of solution accuracy
that is applicable with and without preconditioning, and ensures that
comparisons between different reorthogonalisation strategies in our
experiments are based on a consistent measure of solution quality rather
than on residual-based quantities alone.

\subsection{Reorthogonalisation within GKLB}\label{sec:fullRO}
A detailed analysis of the finite precision behaviour of the symmetric Lanczos process was first provided 
in~\cite{paige1971computation}, primarily in the context of eigenvalue computations. 
This formed the basis for the 
development of a number of effective reorthogonalisation strategies, which we discuss in Section~\ref{sec:part_select}. 
However, GKLB involves two orthogonal bases. In floating-point arithmetic, orthogonality of the columns of both $U_j$ and $V_j$  is typically  lost as $j$ increases. 
In \textit{full reorthogonalisation}, the newly computed vectors $u_{j+1}$ and $v_{j+1}$ 
are reorthogonalised against all previously generated basis vectors (that is, 
$\mathcal{S}^U=\mathcal{S}^V=\{1,2,\ldots ,j\}$  in Algorithm~\ref{alg:LSQR}),
restoring orthogonality to within unit roundoff.  This approach requires 
storing all previous vectors (with {an additional} $j$ vectors of length $n+p$) 
and incurs an additional computational cost 
(approximately $j^2 (n+p)$ flops). Unless the number of LSQR iterations is very small compared to the dimension of the 
problem, this is prohibitively expensive and thus not used in practice.

It is known from the analysis in \cite{paige1971computation}
that loss of orthogonality happens in the direction of the converged Ritz vectors. 
Each step of the GKLB process within LSQR produces a projected bidiagonal 
matrix $B_j$ whose singular values are equal to the Ritz values of $C$.
The computed singular values of $B_j$ determine how quickly different spectral components of the solution are resolved.
Unfortunately, in finite precision arithmetic, loss of orthogonality degrades the quality
of the projected problem. The matrix $B_j$
 may no longer represent an accurate projection of $J$, 
 leading to spurious
duplication of Ritz values. Reorthogonalisation mitigates this behaviour by
maintaining a more accurate projection and, hence, more reliable Ritz values.

It has been shown in \cite{sizh:2000} that only reorthogonalising $v_j$ with respect 
to $V_j$ (and leaving $u_j$ as computed) leads to convergence results 
for LSQR that are comparable to those for full reorthogonalisation, while significantly reducing storage and computation costs.  Our own numerical experiments confirmed this and so, in this study, the default full \RO is
\textit{one-sided reorthogonalisation} with respect to the vectors in $V$.

\subsection{Error bounds for GKLB with \RO}
The effect of \RO  on the quality of the solution returned by LSQR has been  considered recently by Li, Tan and Zhao~\cite{litz:2025}. 
They provide results on the forward-backward error of GKLB with reorthogonalisation of both $U_j$ and $V_j$,  generalising the analysis of Barlow~\cite{barlow:2013} for one-sided GKLB. While the bounds are derived 
without any assumptions on the orthogonalisation strategy, the authors focus on the implications of full reorthogonalisation. 
{Theorem~\ref{th:Li_et_al_3.3} below presents a version of Theorem 3.3 of \cite{litz:2025}, modified for 
one-sided reorthogonalisation with respect to the  $V_j$.
The relationship with the work in~\cite{barlow:2013} is discussed after the proof.} 
The following three quantities are used in the theoretical discussions here and in Section \ref{sec:precision}: 
\begin{itemize}
    \item {true quantities} arising when using a Krylov method to solve problem (\refeq{eq:LSproblem}) in \emph{exact arithmetic};
    \item {perturbed quantities} arising when using a Krylov method to solve the perturbed problem $\min_{\bar{x} \in \mathbb{R}^n} \norm{\bar{b} - \bar{J} \bar{x}}$ in \emph{exact arithmetic} (denoted  using bars);
     \item {computed quantities} arising when using a Krylov method to solve  problem (\refeq{eq:LSproblem}) in \emph{finite precision arithmetic} (denoted using hats).
\end{itemize}

\begin{theorem}
\label{th:Li_et_al_3.3}
Let $J \in \mathbb{R}^{p \times n}$ with $p > n$, and let $b \in \mathbb{R}^p$ be a nonzero
 vector. Suppose $j$ steps of GKLB with one-sided \RO are performed in finite precision
arithmetic with unit roundoff $\bu$. Assume that the rounding errors {in steps~\ref{step:q_for_u}-\ref{step:get_beta} and}
in steps~\ref{step:q_for_v}-\ref{step:get_alpha} 
of Algorithm~\ref{alg:LSQR} are accounted for in vector $f_j$
and matrix $G_{j+1}$, respectively, with {$\norm{f_{i}}= O\!\left(\norm{J}\,\bu\right)$ {for $i\leq j$} and} $\norm{G_{j+1}} = O\!\left(\norm{J}\,\bu\right)$. If the orthogonality level satisfies
$\mu(\bhV_{j+1})<1/2$, then there exist orthonormal matrices $\bar{U}_{j+1}$ and $\bar{V}_{j+1}$,
a perturbation matrix $E^{(j)} \in \mathbb{R}^{p\times n}$ and perturbation vector
$\bub \in \mathbb{R}^p$, such that the following exact GKLB
relations hold:
\begin{gather*}
    \bar{U}_{j+1}\!\left(\widehat\beta_1\, e_1^{(j+1)}\right) = b + \bub, 
    \qquad
    (J + E^{(j)})\,\bar{V}_j = \bar{U}_{j+1}\,\bhB_j, \\
    (J + E^{(j)})^T\bar{U}_{j+1} = \bar{V}_j\,\bhB_j^T
        + \halpha_{j+1}\,\bar{v}_{j+1}\!\left(e_{j+1}^{(j+1)}\right)^{\!T}. \nonumber
\end{gather*}
The perturbations satisfy
\begin{equation} \label{eq:3.24}
    \norm{E^{(j)}} = O\!\left(\norm{J}j\sqrt{n}\,(\bu + \mu(\bhV_{j+1}) ) \right),
   \qquad
    \norm{\bub} = O\!\left(\norm{b}\,\bu\right).
\end{equation}
Furthermore, the computed Lanczos matrix $\bhV_{j+1}$ satisfies the bound
\begin{equation} \label{eq:3.23b}
    \norm{\bhV_{j+1} - \bar{V}_{j+1}} \leq \mu(\bhV_{j+1}) + O\!\left(\mu(\bhV_{j+1})^2\right).   
\end{equation}
\end{theorem}

\begin{proof}
The proof follows the same steps as the proof of Theorem~3.3 of \cite{litz:2025}, with the following changes. Lemma~3.2 of \cite{litz:2025} is replaced by Lemma~\ref{lemma:modified_3.2}. This gives equation (3.12) of \cite{litz:2025} with $\|\widetilde X\| = O\!\left(\|J\|\,j\sqrt{n}\,(\bu+\mu(\bhV_{j+1}))\right)$ and $\mu(\widetilde V)=\mu(\bhV_{j+1})$. The matrix $\bar{U}_{j+1}$ is constructed from the SVD of $\widetilde{P}_{21}$ exactly
as in \cite{litz:2025}. The construction of $\bub$ and $E_1$ in equations
(3.29)--(3.30) uses only $\|[C,\ I_p]\|\le\sqrt2$, and is unaffected. Hence $\|\bub\|= O\!\left(\|b\|\bu\right)$ and $\|E_1\|\le\sqrt2\|\widetilde X\|$. In equation (3.31) of \cite{litz:2025}, $\sigma_{\min}(\widetilde V)\ge(1-2\mu(\bhV_{j+1}))^{1/2}$ follows from $\mu(\widetilde V)=\mu(\bhV_{j+1})<1/2$. This gives~\eqref{eq:3.23b} and
$\|E^{(j)}\|\le\|J\|(\mu(\bhV_{j+1})+O\!\left(\mu(\bhV_{j+1})^2\right)+\sqrt2\|\widetilde X\|$,
which is~\eqref{eq:3.24}. The exact GKLB relations follow by taking the first
$j$ and $j+1$ columns of equation (3.32), as in \cite{litz:2025}.
\end{proof}

The bound \eqref{eq:3.24} depends only on $\mu(\bhV_{j+1})$ at iteration $j$ and requires only a mild assumption on the orthogonality level, {$\mu(\widehat{V}_{j+1})<1/2$}. In contrast, \cite{barlow:2013} uses a reorthogonalisation strategy that requires additional assumptions at each iteration (see  \cite[equations (2.6), (2.7)]{barlow:2013}). The bound in Theorem 5.2 of \cite{barlow:2013} also depends on the orthogonality level of $\bhV_{n}$ in order to bound perturbations at iteration $j$.
 Note that \cite{barlow:2013} considers a GKLB process that returns an upper bidiagonal $B_j$.

Accounting for the finite precision error in the GKLB process and assuming that all other computations in LSQR are performed in exact arithmetic {with $x_0=0$},
 Li et al. 
 \cite{litz:2025} also bound the error between the computed iterate $\hat{x}_j$ and the exact arithmetic iterate $\bar{x}_j$ of the 
perturbed problem, to obtain 
\begin{equation}
    \frac{\norm{\hx_j - \bar{x}_j}}{\norm{\bar{x}_j}} \leq \mu(\bhV_j) + O\!\left(\mu(\bhV_j)^2\right). \label{eq:4.6}
\end{equation}

Theorem~\ref{th:Li_et_al_3.3} shows how the {computed} quantities relate to the perturbed quantities:
the precise relations between the true and perturbed quantities are not established.
Importantly, Theorem~\ref{th:Li_et_al_3.3} holds when, at each step $j$,
the reorthogonalisation is performed against only a subset of previous vectors and regardless of \emph{which} subset, $\mathcal{S}^V$, is used, {provided $\mu(\bhV_{j+1})<1/2$}.
Finally, we note that in the case of full reorthogonalisation, in which 
$\mu(\bhV_j) = O(\bu)$, the perturbation bound (\ref{eq:3.24}) reduces to
$\norm{E^{(j)}}/\norm{J\vphantom{X}} = O(\sqrt{n}\,j\,\bu)$, and the bound in \eqref{eq:3.23b} reduces to $O(\bu)$.
This guarantees that $\hat{x}_j$ is close to $\bar{x}_j$, and the perturbed problem is close to the original problem.

Using the results above, we can now provide a new bound on the relative solution error at iteration $n$.  This requires relating the solution $x_{true}$ of the original least-squares problem and the solution $\bar{x}_{true}=\bar{x}_n$ of the perturbed problem.
Writing \eqref{eq:3.24} with $j=n$ as $\norm{E^{(n)}}\leq \epsilon \norm{J}$  with 
$\epsilon =O\!\left(n^{3/2}\,(\bu + \mu(\bhV_n))\right)$, then if $\epsilon \kappa(J)<1$ and $J$ and $J+E^{(n)}$ are full rank, using Wedin's theorem (see, e.g., 
\cite[Theorem 20.1]{high:02}) and \eqref{eq:4.6} gives
\begin{align}
  \frac{\norm{x_{true} - \hx_n}}{\norm{x_{true}}} = &\, \frac{\norm{x_{true} - \bar{x}_n + \bar{x}_n - \hx_n }}{\norm{x_{true}}} \nonumber \\
                \leq &\, \frac{\norm{x_{true} - \bar{x}_n}}{\norm{x_{true}}} + \frac{\norm{\bar{x}_n - \hx_n }}{\norm{x_{true}}} \nonumber \\
                \leq &\, \frac{\kappa(J) \epsilon}{1 - \kappa(J) \epsilon}\left(2 + (\kappa(J) +1) \frac{\norm{b - J x_{true}}}{\norm{J}\norm{x_{true}}} \right) + \nonumber \\
                 &\, \qquad\left(\mu(\bhV_n) + O\!\left(\mu(\bhV_n)^2\right) \right)\frac{\norm{\bar{x}_n}}{\norm{x_{true}}}. \label{eq:forward_error_bound}
\end{align}
Provided $\norm{\bar{x}_n}/\norm{x_{true}} < \kappa(J)\left(2+\kappa(J)\norm{b - J x_{true}} / \left(\norm{J}\norm{x_{true}}\right)\right)$, this bound is dominated by the first term. 
Thus, 
if $\epsilon = O\!\left(\bu\right)$, the first term is essentially the bound for the forward error of a least-squares problem when using a backward stable solver.

\section{When is employing \RO within LSQR useful?}
\label{sec:Ritz}
Carson, Liesen, and Strako{\v s}~\cite{carson2024} recently presented a detailed comparison of the behaviour of the CG method in exact and finite precision arithmetic for solving systems of linear equations, drawing on a range of computational examples. Their article eloquently summarises the main contributions to this well-studied field and provides an extensive list of references. We highlight one point that is particularly relevant to our current study.
Specifically, Section 2.3 of \cite{carson2024} discusses the sensitivity of CG to
rounding errors. The main finding, based on theoretical results from 
\cite{GREENBAUM1989,GreenbaumStrakos1992,Notay1993}, is the following.
\begin{quote}
{\it Particular eigenvalue distributions, specifically those with
large outlying eigenvalues, cause CG convergence to be more susceptible to
delay caused by finite precision errors.}
\end{quote}

In what follows, we use numerical experiments to illustrate an analogous statement about \RO within LSQR.
\begin{quote}
{\it If the matrix $J$ in (\ref{eq:LSproblem}) has large outlying singular values, the convergence of LSQR is more susceptible to
delay caused by loss of orthogonality in finite precision and using \RO can significantly improve convergence.}
\end{quote}

\subsection{Numerical experiments for LSQR applied to synthetic test cases}\label{sec:SyntheticExamples}
\begin{figure}
    \centering
    \includegraphics{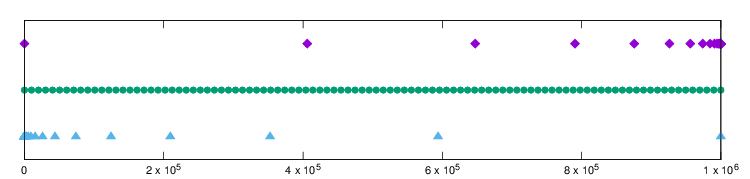} 
    \caption{Singular values corresponding to Case 1 (top line with diamonds), Case 2 (middle line with circles) and Case 3 (bottom line with triangles).}
    \label{fig:CarsonSingularValues}
\end{figure}

We begin by applying LSQR to synthetic examples that are analogous to those in 
\cite{carson2024}, in the sense that the distribution of singular values in our examples is the same as the distribution of eigenvalues in their examples. We choose the dimensions of $J$ to be $p=110$ and $n=100$, and construct $J=Q_p\Sigma Q_n$, where $Q_p\in \mathbb{R}^{p\times p}$ and $Q_n\in\mathbb{R}^{n\times n}$ are random orthogonal matrices. Following \cite[Section 2.1]{carson2024}, we consider three distinct cases defined by the first $n$ diagonal entries of $\Sigma$, with each case having singular values lying in $[0.1,10^6]$. In \textbf{Case 1} a few singular values
are relatively small outliers; in \textbf{Case 2} the singular values are evenly distributed; and in \textbf{Case 
3} a few singular values are relatively large outliers; see Figure~\ref{fig:CarsonSingularValues}.

We define $x_{true} = \frac{1}{\sqrt{n}}[1,1,\ldots, 1]^T$ and generate an error
vector $e$ that is orthogonal to $\operatorname{range}(J)$
     (i.e. $J^Te = 0$), with $\|e\| = 1$. We then construct the vector $b$ in (\ref{eq:LSproblem}) to be $b=Jx_{true}+ \gamma
e$, where $\gamma$ determines the level of noise on the true solution. 
Preliminary experiments found that the effect of reorthogonalisation on 
the performance of LSQR was  qualitatively similar in both a low noise regime ($\gamma = 10^{-6}$) and a high
noise regime ($\gamma = 0.1$), so we present results for the latter only. The norm of the least-squares residual $r = b- Jx$ cannot achieve values smaller than $\gamma$, whereas
the solution error $\|x_{true}-x\|$ and the norm of the normal equation residual $\| J^T r\|$ are
able to achieve relative values much smaller than $\gamma$ (as described in \cite[Section 4]{cptp:2009}).

\begin{figure}[h]
    \centering 
    \includegraphics[
    width=0.94\linewidth,trim = 38mm 0mm 40mm 0mm, clip]{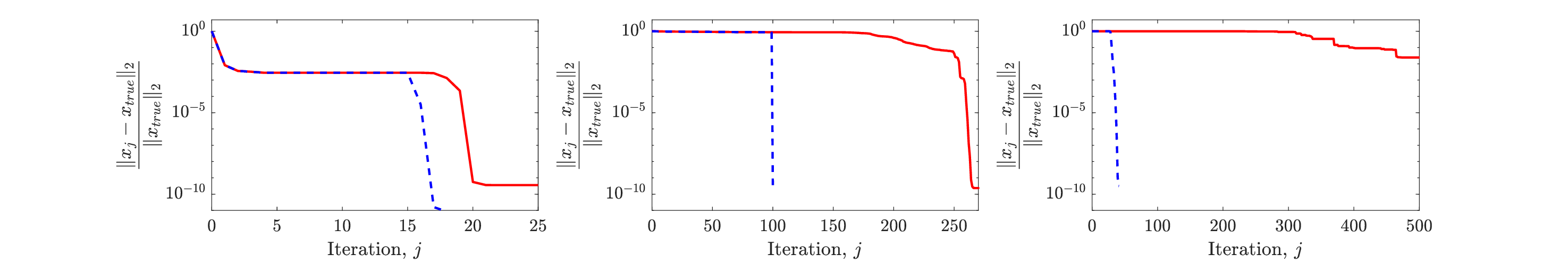}
    \caption{The relative solution error of LSQR for Case 1 (left), Case 2 (middle) and Case 3 (right).
 Each plot shows results with (dashed line)  and without (solid line) full reorthogonalisation.}
    \label{fig:CarsonCases}
\end{figure}

Figure \ref{fig:CarsonCases} shows plots of the logarithm of the relative solution error, which can be readily 
computed for these small examples, against the LSQR iteration index $j$.
The computation  terminates when either (\ref{eq:ratio_PT}) is satisfied with $\delta = 10^{-15}$
or  $\min(\alpha_{j+1},\beta_{j+1}) < 10^{-15}\,\mathrm{estim}(\|J\|$, where $\alpha_{j+1}$ and $\beta_{j+1}$ are defined in Algorithm~\ref{alg:LSQR}.

\begin{figure}
    \centering
    \includegraphics[width=0.94\linewidth,trim = 42mm 0mm 45mm 0mm, clip]{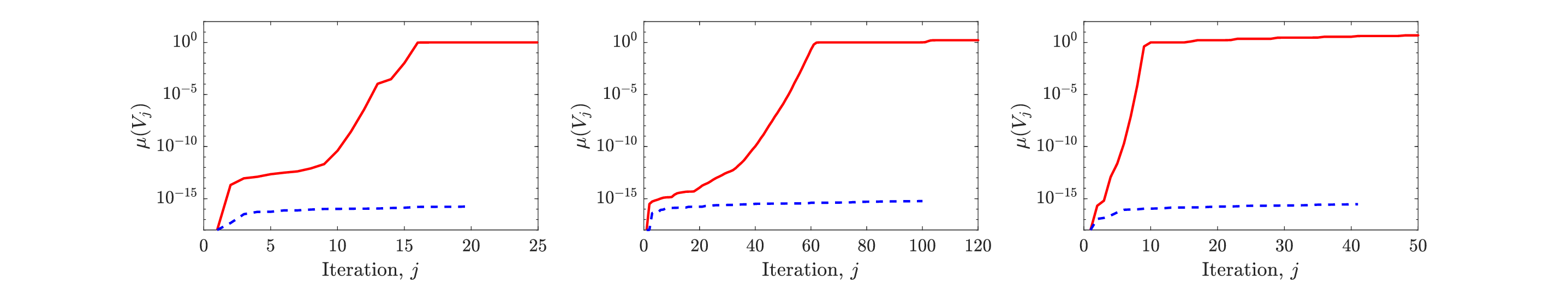}
    \caption{Loss of orthogonality as measured using (\ref{eq:orth index}). 
     Each panel shows results with (dashed line)  and without (solid line) full reorthogonalisation. Note that these figures are zoomed in compared to Figure \ref{fig:CarsonCases}, with the dashed lines ending on termination of LSQR. }
    \label{fig:CarsonMu}
\end{figure}

The general trends in all three cases are analogous to the results in 
\cite{carson2024}. That is, for Case 1,   \RO has little effect; for
Case 2, LSQR requires $n=100$ iterations before the solution error reduces even with reorthogonalisation,
and without reorthogonalisation, there is further stagnation and delay before the requested 
accuracy is achieved; for Case 3, there 
are large discrepancies in the behaviour with and without reorthogonalisation. Specifically, with reorthogonalisation LSQR converges in 41 iterations but, without \RO approximately 2,000 iterations are needed. 
Figure~\ref{fig:CarsonMu} shows plots of $\mu(V_j)$ for each case, which reinforce these observations.

{The orthogonality level $\mu(V_j)$ obtained in the experiments with reorthogonalisation for all Cases satisfy the assumptions of Theorem \ref{th:Li_et_al_3.3}. We recall that the  bound \eqref{eq:forward_error_bound} holds for iteration $n=100$, and can thus only be  used for Case 2.} The right-hand side of \eqref{eq:forward_error_bound}
is dominated by the term $\epsilon\kappa(J)= O\!\left(n^{3/2}\mu(V_n) 10^{7}\right)$. Note that $\mu(V_n)$ stays at $O\!\left(10^{-15}\right)$ (Figure~\ref{fig:CarsonMu}), and thus the final relative error is bounded by $n^{3/2} 10^{-8}\approx 10^{-5}$. While \eqref{eq:forward_error_bound} cannot be applied directly to Case 1 and Case 3, we observe from 
Figure~\ref{fig:CarsonCases} that {in all cases} LSQR returns solutions with relative errors less than $10^{-9}$. 
 
\begin{figure}[hbt]
\centering
\includegraphics[width=0.94\linewidth,trim=20mm 0mm 20mm 0mm,clip]{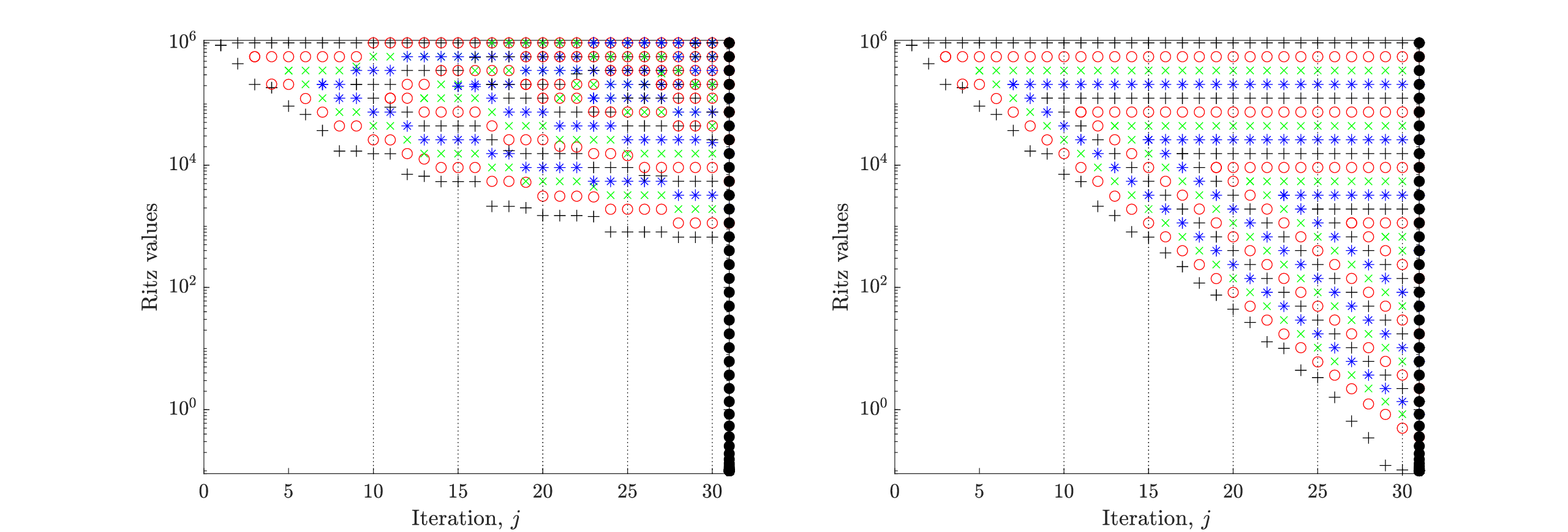}
\caption{Ritz values (plotted on a log scale) for Case 3, without 
(left) and with (right) \RO for the first 30 LSQR iterations. The filled circles represent the singular values of $J$. }
\label{fig:CarsonRitz}
\end{figure}

These experiments suggest that \RO has the most significant effect on LSQR convergence when large outlying singular values 
are present. We now examine {Case 3} in more detail to further support this observation.
Figure~\ref{fig:CarsonRitz} shows the Ritz values without (left) and with (right) reorthogonalisation, 
plotted on a log scale against the iteration index $j$ for the first 30 LSQR iterations; the singular values of $J$ are 
depicted by filled circles on the right-hand side. 
Without reorthogonalisation, the influence of numerical errors on Ritz 
value convergence is clear: as early as iteration 9, Ritz values begin to jump, with multiple values duplicating the largest 
distinct singular values. With reorthogonalisation, by contrast, this duplication is  eliminated, and the Ritz values 
converge smoothly to the individual singular values of $J$. The equivalent plots for {Cases 1} and {2} do not exhibit this 
phenomenon, which is consistent with our earlier observation that \RO offers less benefit in those cases.

{The above experiments show that the presence of large outlying singular values is a key factor in determining how beneficial \RO is for accelerating LSQR convergence.} This has particular relevance in applications where each LSQR iteration is itself costly and only a limited number of iterations are to be performed: in such settings, the use of \RO may yield an overall reduction in total computational time, notwithstanding its added cost per iteration. We now illustrate the potential practical impact of this observation with a test case drawn from variational data assimilation.

\subsection{An illustrative example from data assimilation}

Variational data assimilation (VarDA) is widely used in environmental prediction to ingest observational data into complex and computationally expensive model simulations. It typically involves solving a sequence 
of large-scale sparse linear least-squares problems; see, for example, the recent review article~\cite{dauvzickaite2025introduction}. The operators involved, such as those representing the model simulations, are often too large to be formed or stored explicitly. In numerical weather prediction, these operators are therefore accessed only through matrix--vector products, which account for a major part of the computational cost. Moreover, because of stringent time constraints, the iterative solver must often be terminated after a relatively small number of iterations, before convergence is reached.

In the context of VarDA, loss of orthogonality in CG has been seen in fewer than twenty iterations~\cite{el_akkraoui_preconditioning_2013}. This suggests that finite-precision arithmetic can significantly affect the performance of Krylov solvers, even if the number of iterations is limited. 
In operational weather forecasting, reorthogonalisation is inexpensive compared to matrix--vector products, and is therefore commonly used within CG~\cite{Fisher1998}. \bigskip

\begin{figure}
    \centering
    \includegraphics{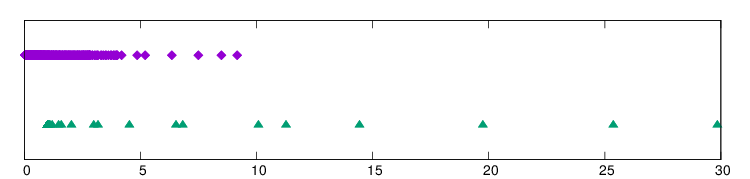}
    \caption{Singular values for {Case 4} corresponding to before (top line with diamonds) and after (bottom line with triangles) preconditioning.} 
    \label{fig:svals_SW}
\end{figure}

\textbf{Case 4} is a VarDA problem based on the one-dimensional shallow water model
of~\cite{lawless2005investigation,lawless2003comparison}, {consisting of two physical variables on a grid of 400 grid points}. {We take a total of 100 observations, 50 of each variable.} After linearisation and standard
reformulation, the resulting system is a least-squares problem 
of the form (\ref{eq:LSproblem}) with $p=900$ and $n=800$ that is solved here using LSQR. 
 The singular values of the associated matrix $J$ are plotted on the top line (with diamonds) in 
Figure~\ref{fig:svals_SW}. 
In VarDA, it is common to solve an equivalent problem with better spectral properties 
by applying  preconditioning {through a variable transformation} (see, e.g.,~\cite[Section 5.1]{dauvzickaite2025introduction} for details). {For this case the condition number of the original problem is small and preconditioning reduces $\kappa(J)$ from 206 to 23. However,  the transformation } 
yields a  pronounced clustering of the singular values of  
{the preconditioned matrix} near unity, which may be more favourable for rapid 
convergence of Krylov subspace solvers, and a number of distinct larger 
singular values (see the lower spectrum in Figure~\ref{fig:svals_SW}). Thus, the preconditioned VarDA problem has a similar 
singular value distribution to Case~3 above.

\begin{figure}[hbt]
\centering
 \includegraphics[width=0.94\linewidth,trim=20mm 0mm 30mm 0mm, clip]{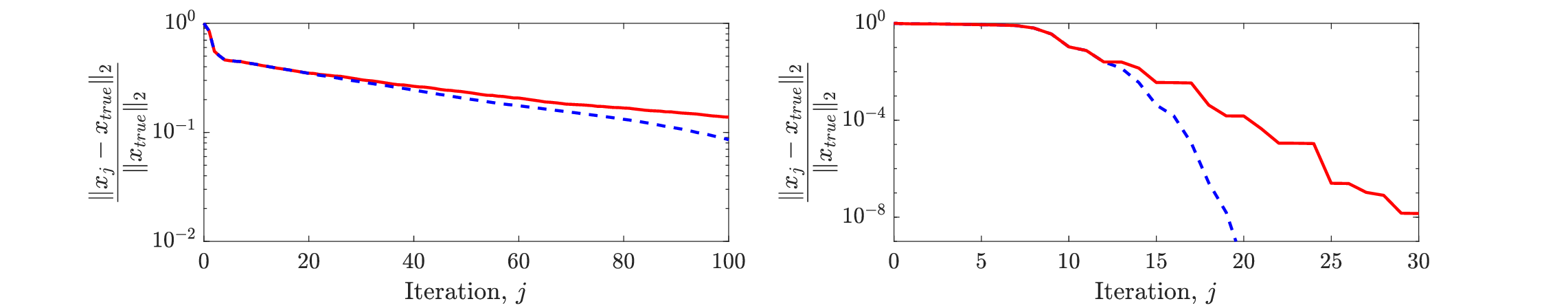}
\caption{
Convergence of LSQR for Case 4  before (left) and after (right) preconditioning has been applied. 
Each plot shows results without (solid line) and with (dashed line) 
full reorthogonalisation.}
\label{fig:SWmotivation}
\end{figure}

Figure~\ref{fig:SWmotivation} shows plots of the logarithm of the solution error {of the original, untransformed variables} against the iteration index for LSQR 
applied to Case 4 before (left) and after (right) preconditioning has been applied (note the different scales on the 
axes).
The acceleration of convergence when preconditioning is employed in
combination with \RO is clear.
The behaviour of Ritz values for the preconditioned case is very similar to that for Case 3 (Figure~\ref{fig:CarsonRitz}) and so plots are not presented here. 

We emphasise that this example highlights an interesting situation where \RO offers no advantage in the unpreconditioned case, but 
is very influential after preconditioning has been applied.

\section{Alternatives to full reorthogonalisation}
\label{sec:local}
So far, we have only used full (one-sided) reorthogonalisation. Other options that are potentially less computationally expensive have been proposed 
including partial and selective reorthogonalisation,
which are described briefly in Section~\ref{sec:part_select}. However, our focus is on large-scale least squares problems for which these methods are impractical
because all the Lanczos vectors must be 
stored. We therefore consider two approaches that require less memory. The first limits
\RO to a small (fixed) number of the previous vectors: this is discussed in 
Section~\ref{sec:ROfixed}. The second utilises mixed precision computations; see Section~\ref{sec:precision}.

\subsection{Partial and selective \RO}
\label{sec:part_select}
It was observed by Paige~\cite{paige1971computation} that full reorthogonalisation of the
Lanczos vectors is unnecessary; instead, maintaining \textit{semi-orthogonality} among the vectors up to 
the square root of unit roundoff is sufficient to guarantee the accuracy of the computed quantities 
in the symmetric Lanczos algorithm up to machine precision~\cite{ParlettScott1979, Simon1984AnalysisOT}. 
Several methods
have been proposed to exploit this property and reduce the computational cost of \RO without significantly compromising accuracy. 
For example, \textit{periodic reorthogonalisation}~\cite{grcar1981, Stewart2002}, performs \RO  against all preceding vectors
whenever the semi-orthogonality condition is violated. In contrast, \textit{partial reorthogonalisation} selects 
only a subset of vectors for reorthogonalisation, based on estimates of the loss of 
orthogonality \cite{simon1984lanczos}.  Furthermore, the analysis in \cite{paige1971computation} on the loss of 
orthogonality in the direction of the converged Ritz vectors led to the idea of
\textit{selective  reorthogonalisation} \cite{ParlettScott1979, Simon1984AnalysisOT}. 
A detailed finite precision study of Lanczos in CG methods in the context of solving linear systems of equations is given 
in~\cite{meurant_lanczos_2006}. In all of these cases,  maintaining semi-orthogonality between the vectors ensures that 
the assumption on the orthogonality level of 
$V$ required by Theorem \ref{th:Li_et_al_3.3} is met.

\subsection{Reorthogonalisation within a fixed memory budget}
\label{sec:ROfixed}
Memory requirements can be fixed \textit{a priori} by only reorthogonalising against subsets of the previously computed 
vectors (the sets $\mathcal{S}^U$ and $\mathcal{S}^V$ in Algorithm~\ref{alg:LSQR}). The 
standard local approach is to reorthogonalise against the most recent $k >0$ vectors. For one-sided reorthogonalisation, 
at iteration $j+1$, $\mathcal {S}^V=\{\max(1,j-k+1),\dots,j\}$  and $\mathcal{S}^U$ is the empty set. This was 
suggested as part of a band Lanczos algorithm~\cite{ruhe1979implementation}. Subsequently,  reorthogonalising the latest block of Lanczos vectors against the two previous blocks was proposed, the key point being that these blocks are still in the core 
memory~\cite{scott1979block}.

A decade ago, extensive tests on practical problems using LSQR {and LSMR} with the standard one-sided local reorthogonalisation 
were reported in~\cite{gosc:2017}.
These experiments observed that a small $k$ can be beneficial (for instance, $k$ in the range 10 to 100 may be an effective choice) but  if a good preconditioner is available, then the advantages of using local \RO were found to be generally less pronounced. 

In Appendix B of a recent manuscript~\cite{chen_preconditioning_2026}, within the context of the conjugate gradient method for symmetric linear systems of equations, the idea of orthogonalising the most recently computed basis vector against the \textit{first} $k$
vectors was proposed, that is, for one-sided reorthogonalisation, at iteration $j+1$, $\mathcal {S}^V=\{1,\dots,\min(k,j)\}$ and $\mathcal{S}^U$ is the empty set. Figure 3.1 of~\cite{simon1984lanczos} indicates that, for the Lanczos algorithm, loss of orthogonality `percolates' backwards, that is, local orthogonality is maintained against the most recent vectors but is lost against earlier ones. This describes the mechanism by which orthogonality is lost, and helps explain why the standard local approach of reorthogonalising against the most recent $k$ vectors can be insufficient to arrest the long-term drift of the basis. This finding has a further theoretical resonance with the analysis of inexact Krylov methods~\cite{Simoncini2003}. In that framework, perturbations introduced at early iterations are the most damaging, since an error at step  $j$ propagates through every subsequent recurrence step, while a late perturbation contaminates far fewer future vectors. Loss of orthogonality against the initial basis vectors $v_1,\ldots,v_{k}$
can be understood as precisely such an early perturbation, whereas loss of orthogonality against recent vectors constitutes a late perturbation with more limited downstream effect. Together, the empirical observation 
in~\cite{simon1984lanczos} and  framework ~\cite{Simoncini2003} provide both a mechanistic and a theoretical rationale for why reorthogonalising against the first $k$ vectors may be more effective than the standard local approach. We note, however, that in~\cite{Simoncini2003}  perturbations to the matrix--vector product in general Krylov methods are analysed, and the mapping to loss of orthogonality in the two-sided GKLB recurrence is by analogy rather than by direct derivation. A rigorous theoretical treatment of this connection in the GKLB setting remains open and would be of independent interest. 

\subsubsection{Numerical experiments with Case 3}
\label{sec:numexpfixed}
\begin{figure}
    \centering
    \includegraphics[width=0.94\linewidth,trim=15mm 0mm 30mm 0mm, clip]{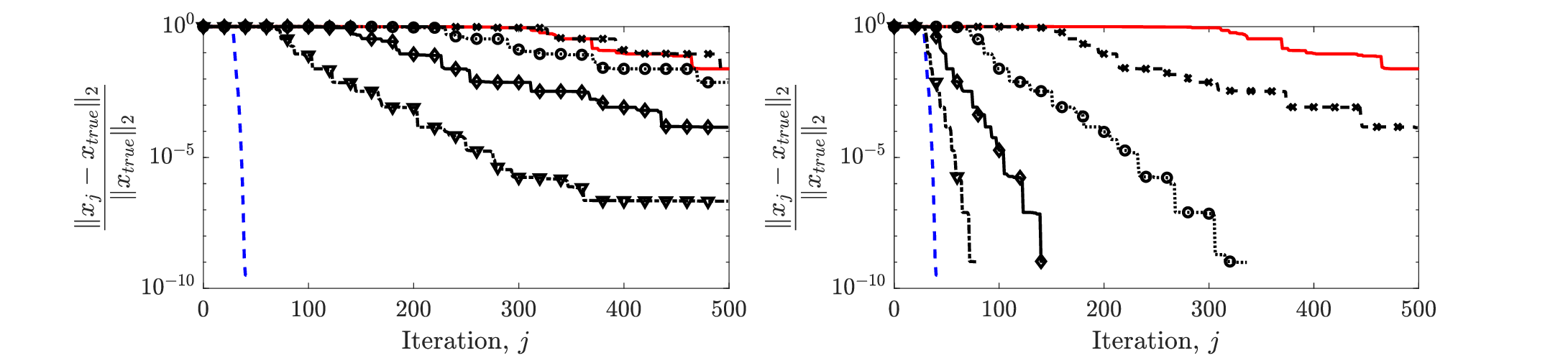}
    \caption{Convergence of LSQR for Case 3. Each plot shows results with (dashed line)  and without (solid)  reorthogonalisation. The left-hand plot also shows results for \RO against the latest $k$ vectors  for $k=5$ (dashed line with crosses), $10$ (dotted line with circles), $15$ (solid line with diamonds) and $20$ (dot-dashed line with triangles). The equivalent results for  \RO against the first $k$ vectors  are shown in the right-hand plot.}
    \label{fig:FixedBudgetConvergence}
\end{figure}

\begin{figure}
    \centering  
   \includegraphics[width=0.94\linewidth,trim=15mm 0mm 30mm 0mm, clip]{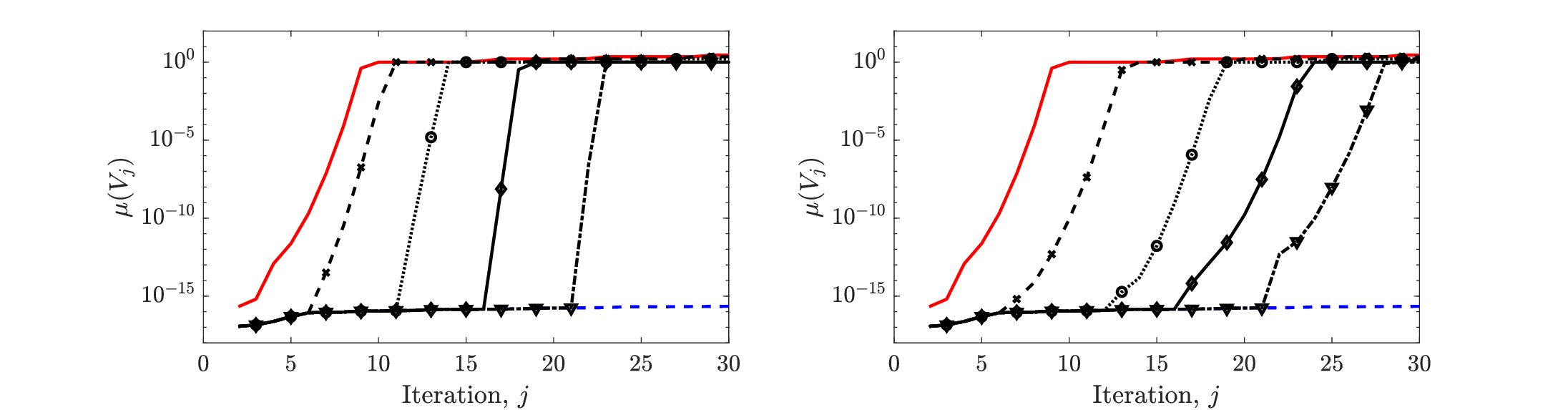}
    \caption{Behaviour of $\mu(V_j)$  for the first $30$ iterations of LSQR for Case 3. Each plot shows results with (dashed line)  and without (solid line)  reorthogonalisation.
    The left plot shows results for  \RO against the latest $k$ vectors for $k=5$ (dashed line with crosses), 
    $10$ (dotted line with circles), $15$ (solid line with diamonds) and $20$ (dot-dashed line with triangles).
    The equivalent results for \RO against the first $k$ vectors  are shown in the right plot. }
    \label{fig:FixedBudgetmu}
\end{figure}

In Section~\ref{sec:SyntheticExamples}, we found that full (one-sided) \RO for the Case~3 test problem gave a clear benefit. 
We now investigate employing
\RO against the most recent $k$ vectors or against the first $k$ vectors. Plots of LSQR convergence for several values of $k$ are presented in
Figure~\ref{fig:FixedBudgetConvergence}. These show that using the first $k$ vectors results in
significantly faster convergence of LSQR compared to using the most recent $k$ vectors.
Figure \ref{fig:FixedBudgetmu} shows the corresponding plots of $\mu({V}_j)$ for the first 30 LSQR iterations. These confirm that 
reorthogonalisation against a fixed number of vectors
delays the increase in $\mu(V_j)$, with the orthogonality level corresponding to reorthogonalisation against the 
first vectors growing more slowly than when reorthogonalising against the latest vectors.
We note that for all choices of $k$ considered here, the orthogonality level exceeds $1/2$ at some particular
iteration, e.g. for $k=10$ at iteration 14 when reorthogonalising against the latest $k$ vectors (see the left panel of
  Figure~\ref{fig:FixedBudgetmu}), and  at iteration 19 when reorthogonalising against the first $k$ vectors (see the right panel of
  Figure~\ref{fig:FixedBudgetmu}).
After this iteration the assumptions on Theorem \ref{th:Li_et_al_3.3} are not satisfied, and hence we cannot apply 
the bound (\ref{eq:forward_error_bound}) at the final iteration. However,  good error reduction can be obtained when 
reorthogonalising against the first $k$ vectors for sufficiently large $k$.

\begin{figure}
    \centering
   \includegraphics[width=0.9\linewidth,trim=20mm 10mm 20mm 0mm,clip]{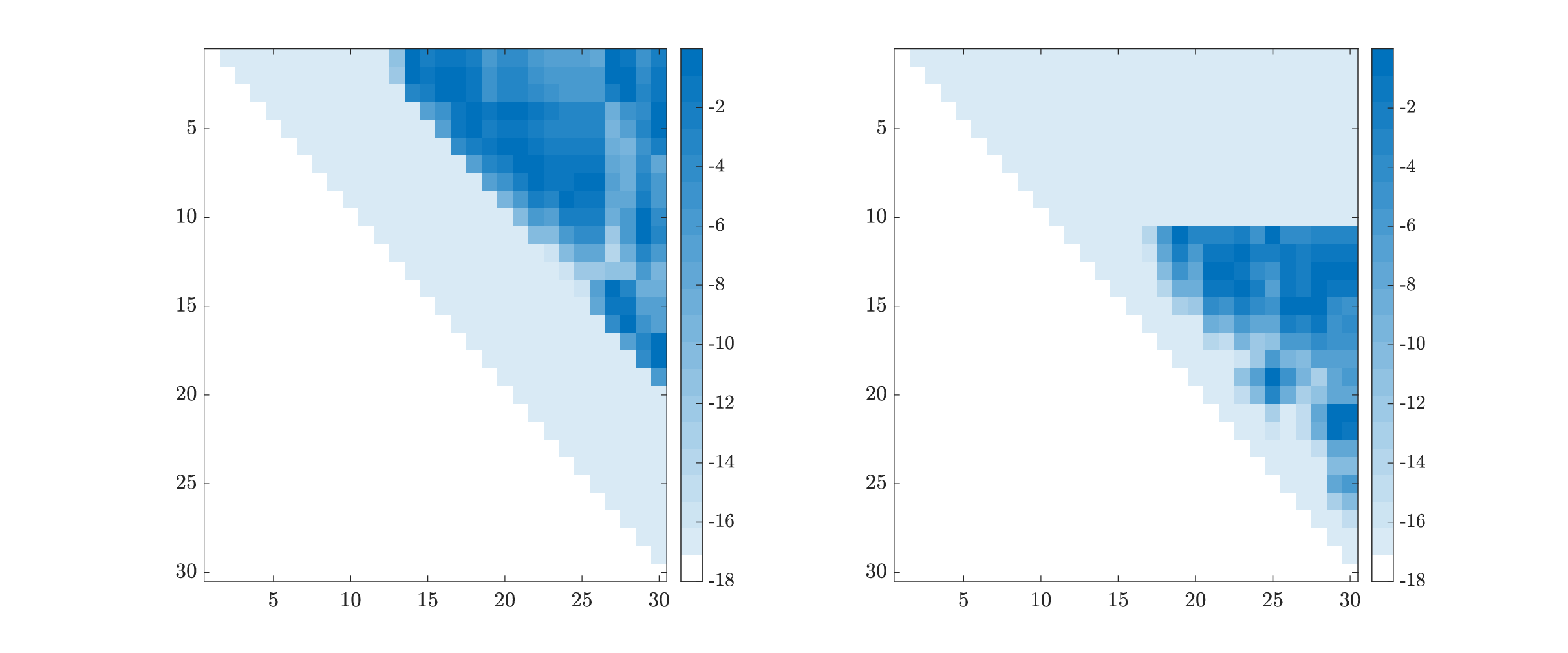}
    \caption{Plot depicting the entries of the matrix $\log |\SUT(I_{30} - V_{30}^TV_{30})|$ with \RO against the latest 10 (left) and 
    first 10 (right) vectors. }
    \label{fig:FixedBudgetTriangle}
\end{figure}

In Figure~\ref{fig:FixedBudgetTriangle}, we 
plot the logarithm of the absolute value of the entries of {the matrix} $\SUT(I_{30} - V_{30}^TV_{30})$ when 
reorthogonalising against the latest $10$ vectors (left panel) and first $10$ vectors (right panel).
We see that \RO against the first $k$ vectors enforces orthogonality to these vectors, while at the 
same time preserving the beneficial `local' orthogonality between adjacent vectors (as seen from the low
values near the diagonal). In this way, the method gains additional `orthogonal' vectors for free. In contrast, \RO against the most recent vectors results in only local orthogonality against the latest $k$ vectors.

\begin{figure}
    \centering
 \includegraphics[width=0.94\linewidth,trim=20mm 0mm 20mm 0mm,clip]{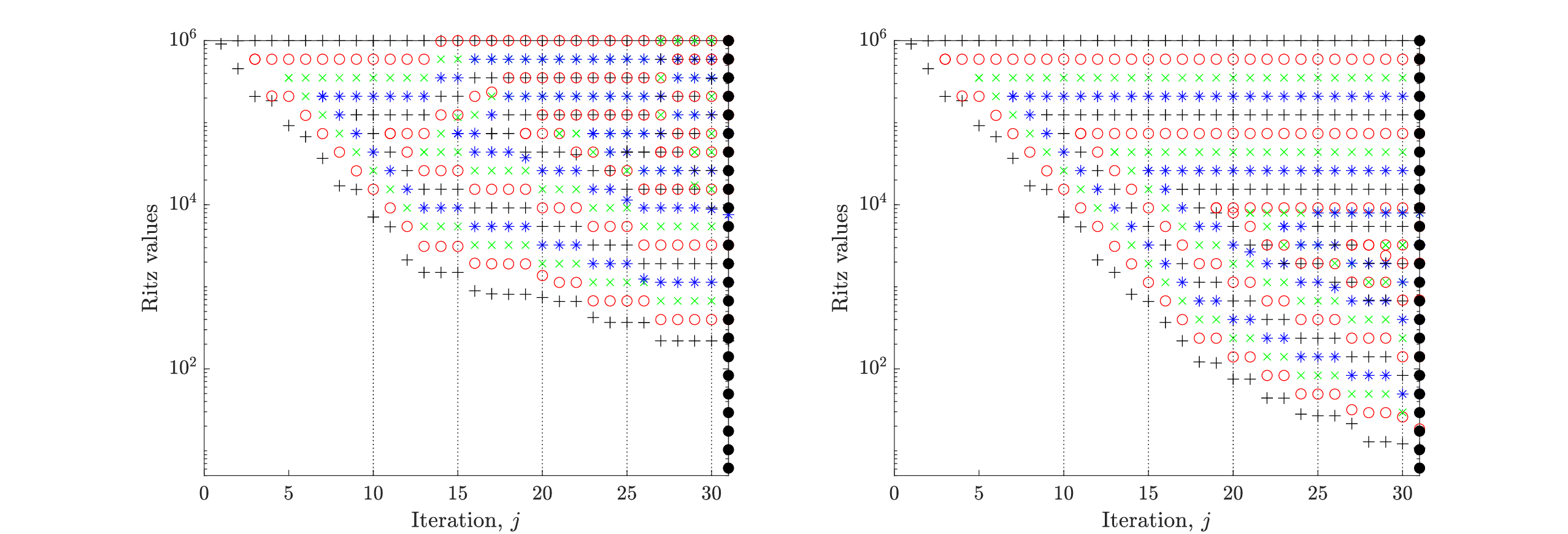}    
    \caption{Ritz values (plotted on a log scale) for Case 3,
with \RO against the latest $10$ vectors (left) and against the first $10$ vectors (right) 
for the first 30 LSQR iterations. 
The solid circles represent the singular values of $J$. }
    \label{fig:FixedBudgetRitz}
\end{figure}

In Figure~\ref{fig:FixedBudgetRitz}, we plot the Ritz values for this case with $k=10$. 
For \RO against the most recent $10$ vectors, the left panel shows that the duplication of Ritz values is delayed compared to when no \RO is used (see Figure \ref{fig:CarsonRitz}). However, after $14$ iterations, the largest Ritz value is duplicated, and this is followed by further duplication at iteration $26$. Duplication of some of the Ritz values also occurs for \RO against the first $10$ vectors (right panel), but it occurs later (iteration 19) and there is no duplication of the largest values (the $10$th largest Ritz value
is the first to be duplicated).

\subsubsection{Numerical experiments with SuiteSparse problems}
\label{sec:suitesparse}
We now consider larger test problems from the study \cite{gosc:2017}, available in the SuiteSparse
Collection\footnote{\url{https://sparse.tamu.edu/}}, using a Fortran implementation of LSQR with the
memory-limited \RO strategies of Section~\ref{sec:ROfixed}.
The \RO uses classical Gram--Schmidt (CGS): $h = V^Tq$, $q \leftarrow q - Vh$, where the columns of $V$ are the stored vectors.
This requires the same work as modified Gram--Schmidt but uses level-2 rather than level-1 BLAS and so
runs much faster. A second pass, needed to guard against cancellation \cite{glr:2005,parl:1998},
is performed only if the first pass reduces $|q|$ by
     more than a factor of $\sqrt{2}$~\cite{dgks:1976}; always performing it
changed the iteration counts by less than 1\% but increased the \RO time by a factor of 1.5--2.4.

The code is compiled with gfortran -O3 and linked to the multithreaded OpenBLAS library; the experiments use
two Intel Xeon E5-2687W processors (16 cores, 20~MB L3 cache per processor). The matrix--vector products are performed serially.
The number of BLAS threads matters: with the default (32), the BLAS threads compete with the
matrix--vector products for cores and memory bandwidth, which slows both down. We use 8 threads, the fastest in our tests.
For ill-conditioned problems the iteration counts can be sensitive to the BLAS implementation and to the number of threads.
Times (in seconds) are measured using {\tt system\_clock}. We use the stopping criterion (\ref{eq:ratio_PT}) with $\delta = 10^{-14}$
(smaller than in practice, to illustrate the effect of \RO\unskip) and $b = \frac{1}{\sqrt{p}}[1,1,\ldots, 1]^T$.

Table~\ref{tab:suitesparse} presents results for mod2 ($66,409 \times 34,774$, {\bf Case~5}), model4 ($4962 \times 1337$, {\bf Case~6}) and
lp\_pilotnov ($2446 \times 975$, {\bf Case~7}). Most of the singular values of $J$ are of order $10-100$, with
outliers of order $10^4-10^7$.
Without \RO\unskip, LSQR fails to converge (Cases 6 and 7) or needs a very large number of iterations (Case 5).
Full \RO gives fast convergence but the memory grows with each iteration. As the required number of iterations is unknown in advance, if convergence is slow this is impractical for large-scale problems.
For both fixed-$k$ strategies the iteration count decreases as $k$ increases, but, as in Section~\ref{sec:numexpfixed},
\RO against the first $k$ vectors is much more effective than against the most recent $k$.
For Case~7, using the most recent $k \le 300$ vectors fails to converge, whereas using the first $k = 300$ needs only 1,674 iterations.
For Case~5, using the most recent $k$ vectors is at best marginally faster than no \RO\unskip, whereas using the first 400 vectors
takes 8.45 seconds, which is faster than full \RO (10.0 seconds), with only 111~MB of storage.
Larger $k$ leads to a slowdown (9.24 seconds for $k = 500$) because the extra \RO cost outweighs the further reduction in the iteration count.

For Case~5, the total time is not monotonic in $k$: with the first $k$ vectors, $k = 50$ and $60$ are faster than no \RO
but $k = 100$ and $200$ are slower. This appears to be a cache effect. The \RO is limited by memory bandwidth, and for $k \le 60$ the stored vectors
(at most 17~MB) fit in the 20~MB L3 cache, whereas for larger $k$ they must be read from main memory on each iteration,
almost doubling the \RO cost per stored vector. Hence $k$ must be large enough for the reduction in the iteration count
to outweigh this extra cost, as it is for $k = 400$.

\begin{table}[htbp]
{
\footnotesize
\begin{center}
\begin{subtable}[t]{.325\textwidth}
\begin{tabular}{|lrrrr|}
\hline
\multicolumn{5}{|c|}{{Case 5}}  \\
\multicolumn{5}{|c|}{{$66,409 \times 34,774$}}  \\

\multicolumn{1}{|c}{$k$}  &
\multicolumn{1}{c}{$iters$} &
\multicolumn{1}{c}{$T_{RO}$} &
\multicolumn{1}{c}{$T_{MV}$} &
\multicolumn{1}{c|}{$T_{tot}$} \\
\hline\Tstrut

0 & 32,632 & 0.00 & 14.9 & 24.4 \\
 $n$ & 1,112 & 9.00 & 0.60 & 9.99 \\
\hline\Tstrut
$50^\ast$ & 31,711 & 8.67 & 15.8 & 35.8 \\
$60^\ast$ & 29,716 & 9.44 & 14.9 & 34.8 \\
$100^\ast$ & 26,838 & 19.6 & 14.2 & 43.3 \\
$200^\ast$ & 21,328 & 42.4 & 11.6 & 61.1 \\
$400^\ast$ & 4,581 & 19.8 & 2.48 & 23.4 \\
\hline\Tstrut
$50^\dag$ & 20,956 & 5.47 & 10.4 & 21.3 \\
$60^\dag$ & 19,772 & 6.04 & 9.91 & 21.0 \\
$100^\dag$ & 17,463 & 14.2 & 9.24 & 28.1 \\
$200^\dag$ & 8,892 & 20.5 & 4.83 & 27.5 \\
$400^\dag$ & 1,734 & 7.11 & 0.94 & 8.45 \\
\hline
\end{tabular}
\end{subtable}
\begin{subtable}[t]{.325\textwidth}
\begin{tabular}{|lrrrr|}
\hline
\multicolumn{5}{|c|}{{Case 6}}  \\
\multicolumn{5}{|c|}{{$4962  \times  1337$}}  \\

\multicolumn{1}{|c}{$k$}  &
\multicolumn{1}{c}{$iters$} &
\multicolumn{1}{c}{$T_{RO}$} &
\multicolumn{1}{c}{$T_{MV}$} &
\multicolumn{1}{c|}{$T_{tot}$} \\
\hline\Tstrut

0 & NC & NC & NC & NC \\
 $n$ & 546 & 0.02 & 0.05 & 0.08 \\
\hline\Tstrut
$50^\ast$ & 37,508 & 0.75 & 3.52 & 11.8 \\
$100^\ast$ & 20,955 & 0.47 & 1.99 & 4.88 \\
$200^\ast$ & 5,053 & 0.16 & 0.48 & 0.81 \\
$300^\ast$ & 1,514 & 0.06 & 0.14 & 0.23 \\
$500^\ast$ & 813 & 0.04 & 0.08 & 0.13 \\
\hline\Tstrut
$50^\dag$ & 16,834 & 0.29 & 1.58 & 3.43 \\
$100^\dag$ & 6,610 & 0.14 & 0.62 & 1.03 \\
$200^\dag$ & 907 & 0.03 & 0.09 & 0.12 \\
$300^\dag$ & 591 & 0.02 & 0.06 & 0.08 \\
$500^\dag$ & 546 & 0.02 & 0.05 & 0.08 \\
\hline
\end{tabular}
\end{subtable}
\begin{subtable}[t]{.325\textwidth}
\begin{tabular}{|lrrrr|}
\hline
\multicolumn{5}{|c|}{{Case 7}}  \\
\multicolumn{5}{|c|}{{$2446 \times 975$}}  \\

\multicolumn{1}{|c}{$k$}  &
\multicolumn{1}{c}{$iters$} &
\multicolumn{1}{c}{$T_{RO}$} &
\multicolumn{1}{c}{$T_{MV}$} &
\multicolumn{1}{c|}{$T_{tot}$} \\
\hline\Tstrut

0 & NC & NC & NC & NC \\
 $n$ & 782 & 0.06 & 0.02 & 0.09 \\
\hline\Tstrut
$100^\ast$ & NC & NC & NC & NC \\
$200^\ast$ & NC & NC & NC & NC \\
$300^\ast$ & NC & NC & NC & NC \\
$400^\ast$ & 4,820 & 0.50 & 0.13 & 0.79 \\
$500^\ast$ & 2,508 & 0.21 & 0.07 & 0.33 \\
\hline\Tstrut
$100^\dag$ & NC & NC & NC & NC \\
$200^\dag$ & 7,100 & 0.30 & 0.19 & 0.78 \\
$300^\dag$ & 1,674 & 0.09 & 0.05 & 0.15 \\
$400^\dag$ & 1,108 & 0.07 & 0.03 & 0.11 \\
$500^\dag$ & 997 & 0.07 & 0.03 & 0.11 \\
\hline
\end{tabular}
\end{subtable}
\end{center}
}
\caption{LSQR results for Cases 5, 6 and 7 with varying amounts of reorthogonalisation. $iters$ is the number of LSQR iterations.
$T_{RO}$, $T_{MV}$ and $T_{tot}$ are, respectively, the times (in seconds)
for the reorthogonalisation, the matrix--vector products with $J$ and $J^T$, and the total solution time.
$^\ast$ and $^\dag$ denote that the most recent and the first $k$ vectors are used, respectively,
in the reorthogonalisation (CGS with the DGKS criterion).
NC denotes convergence not achieved within 50,000 iterations.}
\label{tab:suitesparse}
\end{table}

\section{Mixed precision}\label{sec:precision}
In this section, we consider  reducing the memory requirements
for reorthogonalisation  using mixed precision as an alternative to 
limiting the number of vectors used in reorthogonalisation. 

\subsection{Reorthogonalisation with mixed precision}

For large-scale problems, the storage requirement for the vectors used in reorthogonalisation can become the dominant 
practical constraint. A natural way to alleviate this, without abandoning reorthogonalisation, is to store the Lanczos 
vectors in a lower precision. For example, if single precision (fp32) is used in place of double precision (fp64), the real 
storage for \RO is reduced by a factor of two. 
Li~\cite{li2024double} investigates using mixed precision in LSQR when 
solving discrete ill-posed problems. This work shows that when the noise in the 
right-hand side is large and thus the semi-convergence appears early in the
iteration process, the least-squares solution quality is not affected by using lower precision
for the GKLB with full reorthogonalisation of $U_j$ and $V_j$. Mixed precision approaches for other Krylov
subspace solvers have been proposed in the literature; notably \cite{aliaga2023compressed}
suggests storing the Krylov basis in reduced precision for GMRES.

Here we consider two mixed precision one-sided \RO strategies. 
In our discussion, we denote  the working and lower precision roundoff units by $\bu$ and $\bu_s$, respectively,
with $\bu \ll \bu_s$.
\begin{description}
    \item{MPS} (\textit{mixed precision storage strategy}): 
store the computed vectors $\hat{v}_j$
involved in the \RO in precision $\bu_s$ and perform \RO using the working precision;
    \item{MPSC} (\textit{mixed precision storage+compute strategy}):  store the
$\hat{v}_j$ and perform reorthogonalisation using precision $\bu_s$.
\end{description}

Both strategies reduce the cost of storing the $jn$ \RO vectors $\hat{v}_j$ by at least a factor of two,
while requiring some conversion between the different precisions. 
MPSC potentially offers computational savings on hardware where lower
precision arithmetic is significantly faster than working precision. 

We now present a corollary of Theorem~\ref{th:Li_et_al_3.3} for the MPS strategy.
\begin{corollary}[Mixed precision GKLB]\label{cor:mp_bidiag}
Assume the $j$ steps of GKLB 
satisfy the assumptions of Theorem~\ref{th:Li_et_al_3.3} with the following modification: after each step,
the Lanczos vectors $\hat{v}_i$ are stored in the lower precision giving $\tilde{v}_i$, with
\[
    \norm{\hat{v}_i - \tilde{v}_i} = O(\bu_s),
\]
and used in this form only for the reorthogonalisation corrections
in step~\ref{step:Ro_v} of Algorithm~\ref{alg:LSQR}. 
{We thus assume that the rounding errors
in steps~\ref{step:q_for_v}-\ref{step:get_alpha} 
of Algorithm~\ref{alg:LSQR} are accounted for in matrix $G_{j+1}$ with $\norm{G_{j+1}} = O\!\left(\norm{J}\,\bu_s\right)$
and full one-sided reorthogonalisation achieves $\mu(\hat{V}_{j+1}) = O(\bu_s)$.}
All other operations are performed in working precision $\bu$. Then there exist orthonormal
matrices $\bar{U}_{j+1}$, $\bar{V}_{j+1}$, a perturbation matrix
$E^{(j)} \in \mathbb{R}^{p\times n}$  and a perturbation vector $\bub \in \mathbb{R}^p$ 
such that \eqref{eq:3.24} and \eqref{eq:3.23b}
hold with
\begin{equation}\label{eq:mp_perturbations}
    \norm{E^{(j)}}/\norm{J} = O\!\left(\sqrt{n}\,j\,\bu_s\right), \qquad
    \norm{\bub}/\norm{b} = O\!\left(\bu\right),
\end{equation}
and
\begin{equation}\label{eq:mp_orthog_level}
    \norm{\hat{V}_{j+1} - \bar{V}_{j+1}} = O(\bu_s).
\end{equation}
\end{corollary}

\begin{proof}
{The modified assumptions $\norm{G_{j+1}} = O(\norm{J}\bu_s)$
and $\mu(\hat{V}_{j+1}) = O(\bu_s)$ satisfy all conditions of Theorem~\ref{th:Li_et_al_3.3}}
with $\bu$ replaced by $\bu_s$. The initialisation
$\beta_1 u_1 = b$ is performed in working precision, so
$\norm{\bub} = O(\norm{b}\bu)$ is retained. Substituting into~\eqref{eq:3.24} and
using the fact that $j\bu_s$ dominates $j\bu$ (since $\bu \ll \bu_s$) gives~\eqref{eq:mp_perturbations}. 
The bound~\eqref{eq:mp_orthog_level} follows
from~\eqref{eq:3.23b} with $\mu(\hat{V}_{j+1}) = O(\bu_s)$. 
\end{proof}
 
{We expect that similar
     bounds for MPSC could be obtained} by accounting for the errors in computing inner products in the lower precision $\bu_s$.
As our numerical results in Section~\ref{sec:precisionexp} below illustrate,
in practice MPS is more robust.

\subsection{Numerical examples}
\label{sec:precisionexp}
\begin{figure}
    \centering
    \includegraphics[width=0.95\linewidth,trim = 20mm 5mm 44mm 0mm]{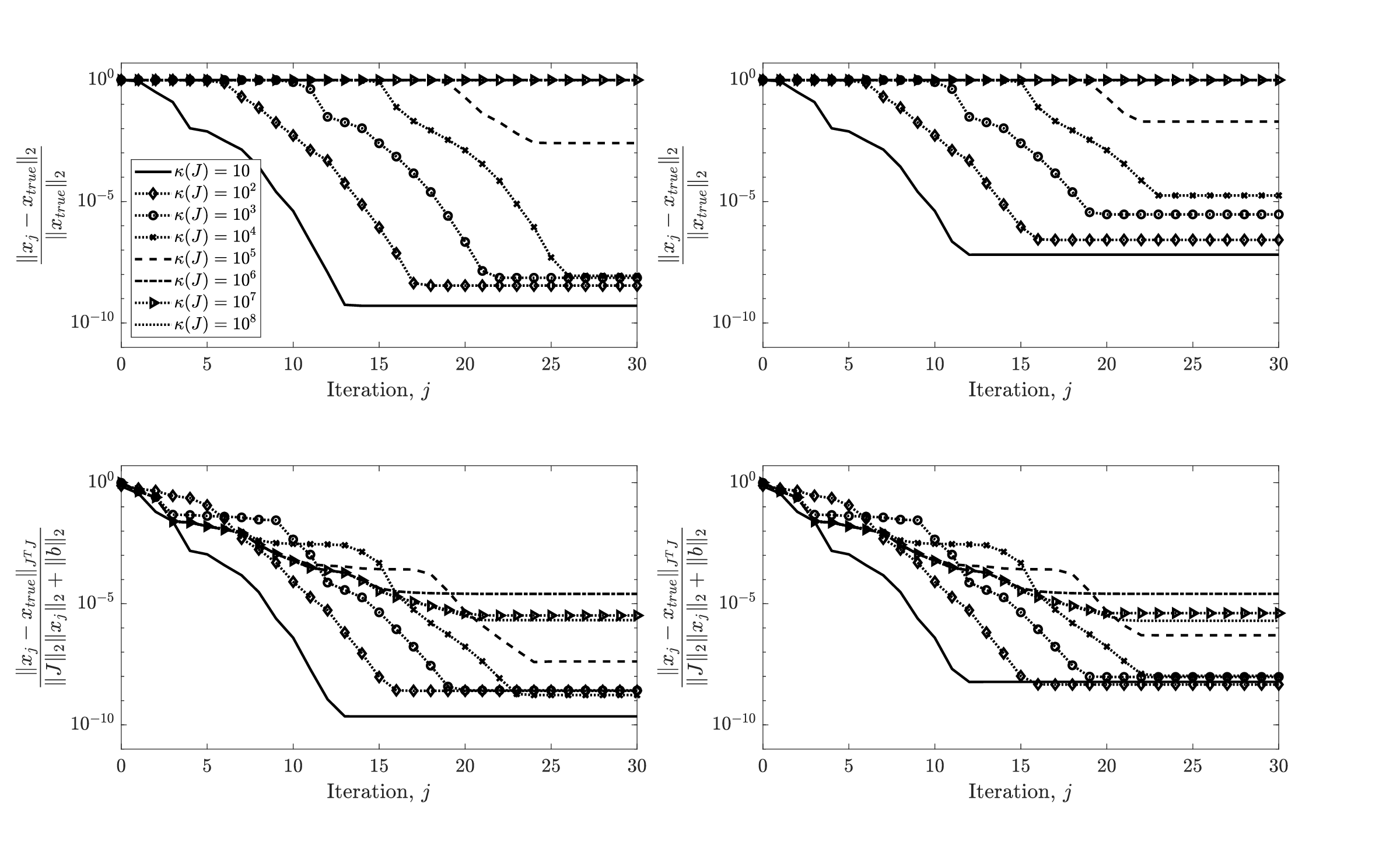}
      \caption{Convergence of LSQR using the MPS  (left plots) and MPSC (right plots) strategies for Case 3 with a range of values of $\kappa(J)$. Full one-sided \RO is used. The forward error (top plots) and the energy norm of the error (bottom plots) are presented. }
    \label{fig:mpRO_forward_error_energy_norm}
    \end{figure}

\begin{figure}
    \centering
        \includegraphics[width=0.95\linewidth,trim = 25mm 5mm 44mm 0mm]{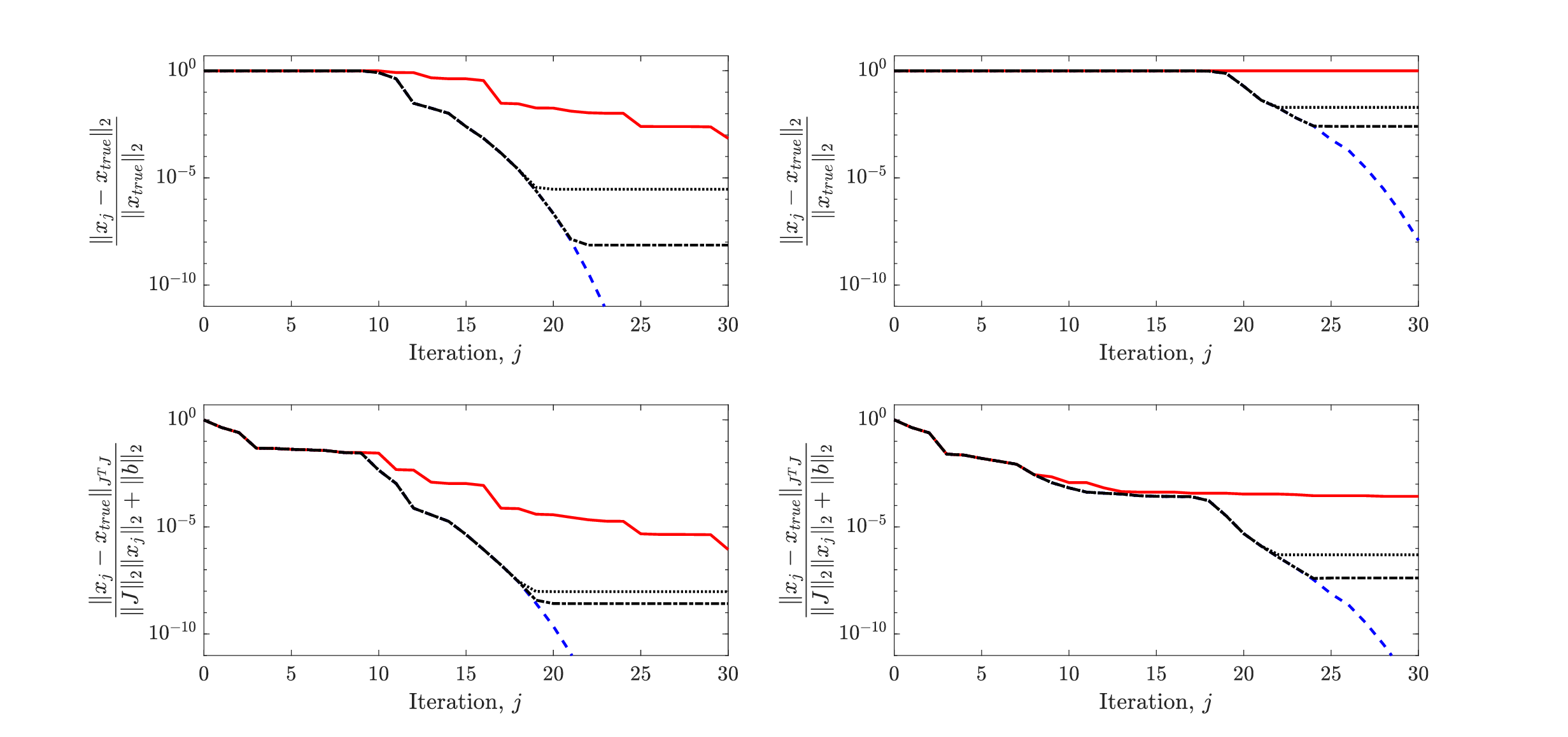}
    \caption{Convergence of LSQR using double precision \RO (dashed), no \RO (solid),  the MPS (dot-dashed) and MPSC (dotted) strategies for Case 3 with $\kappa(J)=10^3$ (left plots) and $\kappa(J)=10^5$ (right plots). 
    The forward error is shown in the top plots and the energy norm of the error is shown in the bottom plots.}
    \label{fig:mpRO_noRO_doubleprecRO}
\end{figure}

\begin{figure}
    \centering
    \includegraphics[width=0.94\linewidth,trim = 24mm 0mm 39mm 0mm,clip]{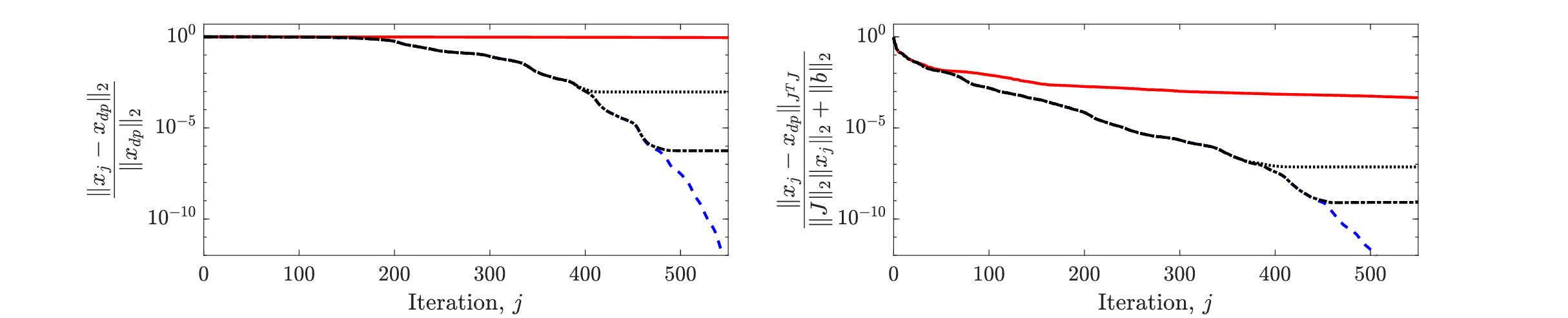}
       \caption{Convergence of LSQR using double precision \RO (dashed), no \RO (solid),  the MPS ({dot-}dashed) and MPSC (dotted) strategies for Case 6 with $\kappa(J)=4.9\times10^4$. 
    $x_{dp}$ is the solution returned by LSQR with \RO in double precision. The forward error is shown in the left plot and the energy norm of the error is shown in the right plot.}
    \label{fig:mpRO_noRO_doubleprecRO_model4}
\end{figure}

We perform numerical experiments using the MPS and MPSC strategies with full (one-sided) reorthogonalisation, with  the working and lower precisions set to
fp64 and fp32, respectively. The least-squares
problems are generated as in Case 3. The largest singular value of $J$ is set
to $1,10^1,\dots,10^7$ so that $\kappa(J)$ ranges from 10 to $10^8$. LSQR convergence curves are presented in
Figure~\ref{fig:mpRO_forward_error_energy_norm}.
These confirm that both strategies are effective over the initial LSQR iterations 
but the error is only improved up to some limiting accuracy that depends on $\kappa(J)$. 
While there is no improvement in the forward error for $\kappa(J)\geq 10^6$, the energy norm of the 
error reduces to around $10^{-5}$, even for highly ill-conditioned problems. Note that MPS
tends to reach a lower limiting accuracy compared to MPSC; the MPS limiting accuracy is similar for
the test cases with $\kappa(J)<\bu_s^{-1/2}$. 
Figure~\ref{fig:mpRO_noRO_doubleprecRO} shows that
initially both 
strategies have the same effect on LSQR convergence as double precision
reorthogonalisation;  differences only appear
when the limiting accuracy of the solution is approached.
In Figure~\ref{fig:mpRO_noRO_doubleprecRO_model4}, similar behaviour is seen 
for Case 6 (problem model4 from SuiteSparse). 

In these experiments, LSQR was run for 30 iterations. 
If the stopping criterion \eqref{eq:ratio_PT} is employed, then
 recall that it can be trusted until the solution reaches its maximal attainable accuracy.
When using the MPS and MPSC strategies, we observe that the ratio \eqref{eq:ratio_PT} continues to decrease after the solution has reached
its attainable accuracy. The parameter $\delta$ should thus  be set large enough to terminate LSQR before the forward error norm stagnates. 

\section{Concluding remarks}
In this paper, we have revisited the role of reorthogonalisation in LSQR from the perspective of 
modern large-scale least-squares applications. Our numerical experiments show that
reorthogonalisation can substantially accelerate convergence when the coefficient matrix has large 
outlying singular values. We have also adapted to LSQR two approaches for reducing the 
cost of reorthogonalisation: a memory-limited approach based on reorthogonalising against the 
first $k$ Lanczos vectors, which consistently outperforms the standard method of using the most 
recent vectors, and two mixed-precision storage and computation strategies that significantly reduce 
memory requirements while preserving the convergence benefits of reorthogonalisation. {We have also presented bounds for the GKLB process and the LSQR solution error
  that apply to any one-sided reorthogonalisation strategy, and depend only on the
  level of orthogonality that the strategy maintains, not on which vectors are used.
  They therefore cover memory-limited and mixed precision reorthogonalisation while $\mu(\hat{V}_{j+1})<1/2$. In our experiments, reorthogonalising against a fixed
  number of vectors eventually violates this condition, so explaining the observed
  convergence of these strategies requires further analysis. 

Our results 
suggest that reorthogonalisation should be viewed not only as a mechanism for maintaining 
numerical stability in finite precision arithmetic, but also as a practical tool for improving the 
efficiency of LSQR in applications where matrix–vector products are expensive and only a limited 
number of iterations can be afforded. 

We  have concentrated on LSQR. The more recent 
LSMR method of Fong and Saunders~\cite{fosa:2011} is  mathematically equivalent to MINRES applied to the normal equations. Because LSMR is
based on GKLB, it faces the same issues of loss of orthogonality as LSQR.
Experiments illustrating this are reported in \cite{fosa:2011, gosc:2017}.

In our study, we have investigated the role of reorthogonalisation in LSQR mainly based on the singular value distribution. As the convergence behaviour of LSQR also depends on the initial residual, a more comprehensive analysis that incorporates this factor is a possible direction for future work.
Other possible avenues to explore include the investigation of adaptive precision techniques in which the Lanczos vectors are stored in varying precisions to reduce memory requirements while preserving the benefits of \RO in working precision \cite{Higham_Mary_2022}. Additionally, the combination of mixed precision and partial \RO strategies could be examined. Inspired by recent work on randomised orthogonalisation for Krylov methods, LSQR could employ low-dimensional sketches of the Golub–Kahan-Lanczos basis to perform \RO approximately in a sketched space, reducing computational and communication costs while maintaining numerical stability with high probability \cite{damas2025a}. 

Finally, we note that there are other 
iterative methods that can be used to solve least-squares problems. For example, \cite{hayi:2010} states that AB-GMRES should be
numerically more stable than CGLS and LSQR 
because it is based on the Arnoldi process, rather than 
on short-term recurrences \cite{hahm:2022}. The authors conjecture that AB-GMRES may be numerically equivalent to CGLS and LSQR with full reorthogonalisation. It would be of interest to explore this further.

\bigskip

\noindent\textbf{Acknowledgements.}
This work was partially funded by a Network Support
Grant (Isaac Newton Institute for the Mathematical Sciences and the Engineering
\& Physical Sciences Research Council (EPSRC), EP/V521929/1) and a Research In
Groups Grant (International Centre for Mathematical Sciences) in the UK.
We also acknowledge funding provided by the Deutsche
Forschungsgemeinschaft (DFG) - Project-ID 318763901 - SFB1294 to visit the
University of Potsdam. 

For the purpose of open access, the authors have
applied a Creative Commons Attribution (CC BY) licence to any Author Accepted
Manuscript version arising.

\bigskip
\noindent\textbf{Author Contributions} All authors contributed equally to the conception, analysis and writing of this manuscript. All authors read and approved the final manuscript.

\bigskip
\noindent\textbf{Data availability} The sparse matrices reported on in this study are included in
the SuiteSparse Matrix Collection \url{https://sparse.tamu.edu/}

\section*{Declarations}
\noindent\textbf{Conflict of interest.} All authors declare that they have no conflict of interest.

\appendix
\section{Lemmas for the proof of Theorem~\ref{th:Li_et_al_3.3}}\label{app:proof_one_sided_ro}

Using no reorthogonalisation of the $\hu_j$ vectors corresponds to setting the reorthogonalisation coefficients $\xi_{ik}$ in equation (2.3) of~\cite{litz:2025} to zero.
In the proof, we use the following relations that account for this
(these are modifications of equations (2.3) and (2.19), respectively, in \cite{litz:2025}):
\begin{align}
    J \hv_j &= \halpha_j \hu_j + \hbeta_{j+1}\hu_{j+1} + f_j,  \label{eq:Jvi} \\
    \bar{c}_j & = \begin{bmatrix}
        \xi_{1j} & \dots & \xi_{jj}
    \end{bmatrix}^T = \  0. \label{eq:cbar=0} 
\end{align}

We state a modified version of Lemma 3.1 of \cite{litz:2025} that holds for one-sided GKLB. The notation has been adapted to that of this paper.
\begin{lemma}[Modified Lemma 3.1 of \cite{litz:2025}]\label{lemma:modified_3.1}
Consider the $j$th step of GKLB with one-sided reorthogonalisation of $\hv_i$. Let
\begin{gather*}
    \check{P}_{j+1} = P_1 \dots P_{j+1}
\end{gather*}
be a product of Householder transformations $P_i = I_{p+n+1} - p_ip_i^T$, where $p_i = \begin{bmatrix}
    -e_i^{(n+1)} \\ \hu_i
\end{bmatrix}$. Then the following holds:
\begin{equation*}
    \check{P}_{j+1}  \begin{bmatrix}
\halpha_j e_j^{(j)}\\
\hbeta_{j+1}e_1^{(p+n-j+1)}
\end{bmatrix} = \begin{bmatrix}
0_{n+1}\\
J \hv_j
\end{bmatrix} + x_j,
\end{equation*}
where 
\begin{equation}
    \Vert x_j \Vert = O\!\left( j\Vert J \Vert  (\bu +\mu(\bhV_j) ) \right).
\end{equation}
\end{lemma}
\begin{proof}
    As in \cite{litz:2025}, we start by computing $ P_{j+1}  \begin{bmatrix}
\halpha_j e_j^{(j)}\\
\hbeta_{j+1}e_1^{(p+n-j+1)}
\end{bmatrix}$. Using \eqref{eq:cbar=0}, we have 
\begin{equation*}
    P_{j+1}  \begin{bmatrix}
\halpha_j e_j^{(j)}\\
\hbeta_{j+1}e_1^{(p+n-j+1)}
\end{bmatrix} = \begin{bmatrix}
\halpha_j e_j^{(n+1)}\\
J \hv_j - \halpha_j \hu_j
\end{bmatrix}
-
\begin{bmatrix}
0_{n+1}\\
f_j
\end{bmatrix}
\end{equation*}
and
\begin{equation}
P_j
\begin{bmatrix}
\halpha_j e_j^{(n+1)}\\
J \hv_j - \halpha_j \hu_j
\end{bmatrix}
=
\begin{bmatrix}
0_{n+1}\\
J \hv_j
\end{bmatrix}
-
\left(
\hbeta_{j+1}\hu_j^{T}\hu_{j+1}
+
\hu_j^{T}f_j
\right) p_j.
\end{equation}
Then 
\begin{align*}
    P_j P_{j+1} \begin{bmatrix}
\halpha_j e_j^{(j)}\\
\hbeta_{j+1}e_1^{(p+n-j+1)}
\end{bmatrix} & \ = \begin{bmatrix}
0_{n+1}\\
J \hv_j
\end{bmatrix} - \left(
\hbeta_{j+1}\hu_j^{T}\hu_{j+1}
+
\hu_j^{T}f_j
\right) p_j
- \begin{bmatrix}
0_{n+1}\\
f_j
\end{bmatrix} + \hu_j^{T}f_j p_j \\
& \ = \begin{bmatrix}
0_{n+1}\\
J \hv_j
\end{bmatrix} - \left(
\hbeta_{j+1}\hu_j^{T}\hu_{j+1}
\right) p_j
- \begin{bmatrix}
0_{n+1}\\
f_j
\end{bmatrix} .
\end{align*}
We now have to obtain a bound on the norm of 
\begin{equation}\label{eq:omega_def}
    \omega_j = - \left(
\hbeta_{j+1}\hu_j^{T}\hu_{j+1}
\right) p_j
- \begin{bmatrix}
0_{n+1}\\
f_j
\end{bmatrix}.
\end{equation}
 We first consider the scalar $\hbeta_{j+1}\hu_j^{T}\hu_{j+1}$. Rearranging \eqref{eq:Jvi} gives \begin{equation*}
    \hbeta_{j+1}\hu_{j+1} = J \hv_j - \halpha_j \hu_j - f_j.
\end{equation*}
Plugging the above into $\hbeta_{j+1}\hu_j^{T}\hu_{j+1}$, and recalling that $\| \hu_j \| =1$, yields
\begin{equation*}
\hbeta_{j+1} \hu_j^T \hu_{j+1} = \hu_j^T(\hbeta_{j+1} \hu_{j+1})   
                          = \hu_j^T(J \hv_j -\halpha_j \hu_j - f_j)
                          = \hu_j^T J \hv_j - \halpha_j - \hu_j^T f_j.
\end{equation*}
Using $J^T \hu_j = \bhV_j \widetilde{d}_j + g_j$ (equation (2.16) in \cite{litz:2025}), where {$g_j$ is the $j$th column of $G_{j}$,} $\widetilde d_j
=
\begin{bmatrix}
\eta_{1j} & \dots & \eta_{j-2,j} & \hbeta_j+\eta_{j-1,j} & \halpha_j
\end{bmatrix}^T$ and  $\eta_{ij}$ are coefficients arising in the reorthogonalisation of $\hv_j$, we obtain
\begin{align}
\hbeta_{j+1} \hu_j^T \hu_{j+1} & \ = (\bhV_j \widetilde{d}_j + g_j)^T \hv_j - \halpha_j - \hu_j^T f_j \nonumber \\
                        & \ = \widetilde{d}_j^T \bhV_j^T \hv_j  + g_j^T \hv_j - \halpha_j - \hu_j^T f_j \label{eq:beta_uT_u}.
\end{align}
We now examine the term $\widetilde{d}_j^T \bhV_j^T \hv_j$.
Notice that
\begin{equation*}
\bhV_j^T \hv_j
=
\begin{bmatrix}
\hv_1^T\\
\vdots\\
\hv_{j-1}^T \\
\hv_j^T
\end{bmatrix}\hv_j
=
\begin{bmatrix}
\hv_1^T \hv_j\\
\vdots\\
\hv_{j-1}^T \hv_j\\
1
\end{bmatrix}.
\end{equation*}
We can thus write
\begin{equation}\label{eq:djVTv}
    \widetilde{d}_j^T \bhV_j^T \hv_j = \halpha_j + \sum_{i=1}^{j-1} [\widetilde{d}_j]_i \hv_i^T \hv_j, 
\end{equation}
where $[\widetilde{d}_j]_i$ denotes the $i$th element of $\widetilde{d}_j$. We now have to bound the scalar $\sum_{i=1}^{j-1} [\widetilde{d}_j]_i \hv_i^T \hv_j$.
{First, recall that $-\begin{bmatrix}
    \hv_1^T \hv_j & \dots & \hv_{j-1}^T \hv_j
\end{bmatrix}^T$ is the $j$th column of $\operatorname{SUT}
\left(I_j-\bhV_j^T \bhV_j\right)$ and thus $\left \| \begin{bmatrix}
    \hv_1^T \hv_j & \dots & \hv_{j-1}^T \hv_j
\end{bmatrix}^T \right \| \leq \mu(\bhV_j)$. Then using Cauchy-Schwarz inequality we have
\begin{equation*}
\left| \sum_{i=1}^{j-1}[\widetilde d_j]_i\,\hv_i^T \hv_j \right|
\leq \| \widetilde d_j \| \mu(\bhV_j).
\end{equation*}
}

Combining this with $\|\widetilde d_j\| \leq \| J \| + O(\| J \| (\bu +\mu(\bhV_j))$ (equation (2.17) of \cite{litz:2025}), we obtain
\begin{equation}\label{eq:bound_sum_dvivj}
\left| \sum_{i=1}^{j-1} [\widetilde d_j]_i\,\hv_i^T\hv_j \right| \le  \mu(\bhV_j)
\left( \|J\|+O(\|J\|(\bu+ \mu(\bhV_j)) \right).
\end{equation}
Coming back to \eqref{eq:beta_uT_u} and using \eqref{eq:djVTv} and \eqref{eq:bound_sum_dvivj}, we obtain
\begin{align*}
\hbeta_{j+1} \hu_j^T \hu_{j+1} 
& \ =
\halpha_j + \sum_{i=1}^{j-1} [\widetilde d_j]_i \hv_i^T \hv_j +  g_j^T \hv_j - \halpha_j - \hu_j^Tf_j\\
&\ = \sum_{i=1}^{j-1} [\widetilde d_j]_i \hv_i^T \hv_j + g_j^T\hv_j - \hu_j^T f_j.
\end{align*}
Using $\| p_j \| = \sqrt{2}$, $\| \hu_i \|, \| \hv_i \| =1$, the above, \eqref{eq:bound_sum_dvivj} and $\|f_j\|,\|g_j\|=O(\|J\|\bu)$, we can obtain a bound for the norm of $\omega_j$ in \eqref{eq:omega_def} as follows:
\begin{align*}
 \| \omega_j\| 
&\ \le \left|\hbeta_{j+1} \hu_j^T\hu_{j+1}\right|\sqrt{2} +\|f_j\|\\
&\ \le \left( \left| \sum_{i=1}^{j-1} [\widetilde d_j]_i \hv_i^T \hv_j \right|
+\|g_j\|+\|f_j\| \right)\sqrt{2} + \|f_j\|\\
&\ \le \sqrt{2} \mu(\bhV_j) \left( \|J\|+O(\|J\|(\bu+\mu(\bhV_j))) \right)
+ O(\|J\|\bu)(2\sqrt{2}+1)\\
& \ = O\!\left( \|J\|(\bu+\mu(\bhV_j))\right).   
\end{align*}
The rest of the proof is as in \cite{litz:2025}, except using the above bound in the final inequality.
\end{proof}

The modification to Lemma 3.2 of \cite{litz:2025} concerns the bound for $\| \widetilde{X} \|$. 

{
\begin{lemma}[Modified Lemma 3.2 of \cite{litz:2025}]
\label{lemma:modified_3.2}
Consider the $j$th step of GKLB with one-sided reorthogonalisation of $\hv_j$. There exist unit vectors $\widetilde{u}_{j+2}, \dots, \widetilde{u}_{n+1}$ and orthonormal vectors $\widetilde{v}_{j+2}, \dots, \widetilde{v}_n$, 
such that $\widetilde{v}_i$ are orthogonal to $\bhV_{j+1}$. Define $\widetilde{U}=\begin{bmatrix}
    \bhU_{j+1} & \widetilde{u}_{j+2} & \dots & \widetilde{u}_{n+1}
\end{bmatrix}$, $\widetilde{V}=\begin{bmatrix}
    \bhV_{j+1} & \widetilde{v}_{j+2} & \dots & \widetilde{v}_n
\end{bmatrix}$ and 
$\widetilde{B}$ to be an $(n+1) \times n$ bidiagonal matrix with $\halpha_1,\dots,\halpha_{j+1},\widetilde{\alpha}_{j+2}, \dots, \widetilde{\alpha}_n$ on the main diagonal and $\hbeta_2,\dots,\hbeta_{j+1},\widetilde{\beta}_{j+2}, \dots, \widetilde{\beta}_{n+1}$ on the subdiagonal, where $\widetilde{\alpha}_i$ and $\widetilde{\beta}_i$ are some positive numbers. Then there exists a matrix $\widetilde{X} \in \mathbb{R}^{p+n+1\times n}$ such that
\begin{equation*}
    \begin{bmatrix}
0_{(n+1) \times n}\\
J \widetilde{V}
\end{bmatrix} + \widetilde{X} = 
\widetilde{P} \begin{bmatrix}
\widetilde{B}\\
0_{p \times n}
\end{bmatrix},
\end{equation*}
where
\begin{equation*}
    \| \widetilde{X} \| = O \left(j \sqrt{n} \| J \| (\bu + \mu(\bhV_{j+1})) \right)
\end{equation*} 
$ \widetilde{P} = \check{P}_{j+1}  \widetilde{P}_{j+2} \dots \widetilde{P}_{n+1}$,
$\widetilde{P}_i = I_{p+n+1} - \widetilde{p}_i \widetilde{p}_i^T$, and $\widetilde{p}_i = \begin{bmatrix}
    -e_i^{(n+1)} \\ \widetilde{u}_i
\end{bmatrix}$. 

Also, let $\widetilde{P} = \begin{bmatrix}
    \widetilde{P}_{11} & \widetilde{P}_{12}\\ \widetilde{P}_{21} & \widetilde{P}_{22}
\end{bmatrix}$, then 
$\widetilde{P}_{21}  = \widetilde{U}(I_{n+1} - \widetilde{P}_{11})$.
\end{lemma}

\begin{proof}
The proof follows the same steps as the proof of Lemma 3.2 of \cite{litz:2025}, with the following changes. The orthogonalisation against $\bhU_{j+1}$ is dropped in equation (3.17) of \cite{litz:2025},
and we set $\widetilde{\beta}_{i+1}\widetilde{u}_{i+1}=J\widetilde{v}_i-\widetilde{\alpha}_i\widetilde{u}_i$ exactly for $i \ge j+1$. The construction of $\widetilde{v}_{i+1}$ is kept, so that $\mu(\widetilde V)=\mu(\bhV_{j+1})$.
The columns of $\widetilde{X}$ are bounded by Lemma~\ref{lemma:modified_3.1} and by the relations in its proof, which involve only $\mu(\bhV_{j+1})$ and not $\mu(\bhU_{j+1})$. Theorem~3.2 of \cite{litz:2025} holds with the new construction of $\widetilde{U}$ and so $\widetilde{P}_{21}=\widetilde{U}(I-\widetilde{P}_{11})$ holds.
\end{proof}
}
\bibliography{DA}
\end{document}